\documentclass{amsproc}
\usepackage{color}
\newtheorem{Theorem}{Theorem}
\newtheorem{theorem}{Theorem}[section]
\newtheorem{lemma}[theorem]{Lemma}
\newtheorem{proposition}[theorem]{Proposition}
\newtheorem{corollary}[theorem]{Corollary}

\theoremstyle{definition}
\newtheorem{definition}[theorem]{Definition}

\theoremstyle{remark}

\numberwithin{equation}{section}

\newcommand{\abs}[1]{\lvert#1\rvert}
\newcommand{\Abs}[1]{\lVert#1\rVert}

\newcommand{\bd}{ {\sc Proof}.\ \ }

\begin{document}

\title{Isolated circular  orders of $PSL(2,\Z)$}


\author{Shigenori Matsumoto}
\address{Department of Mathematics, College of
Science and Technology, Nihon University, 1-8-14 Kanda-Surugadai,
Chiyoda-ku, Tokyo, 101-8308 Japan}
\email{matsumo@math.cst.nihon-u.ac.jp}

\thanks{The author is partially supported by {{Grant-in-Aid for
Scientific Research}} (C) No.\ 25400096.}

\date{\today}

\newcommand{\AAA}{{\mathcal A}}
\newcommand{\BBB}{{\mathbb B}}
\newcommand{\LL}{{\mathcal L}}
\newcommand{\MCG}{{\rm MCG}}
\newcommand{\PSL}{{\rm PSL}}
\newcommand{\R}{{\mathbb R}}
\newcommand{\Z}{{\mathbb Z}}
\newcommand{\XX}{{\mathcal X}}
\newcommand{\per}{{\rm per}}
\newcommand{\N}{{\mathbb N}}

\newcommand{\PP}{{\mathcal P}}
\newcommand{\RR}{{\mathcal R}}
\newcommand{\GG}{{\mathbb G}}
\newcommand{\FF}{{\mathcal F}}
\newcommand{\EE}{{\mathbb E}}
\newcommand{\BB}{{\mathbb B}}
\newcommand{\CC}{{\mathcal C}}
\newcommand{\HH}{{\mathcal H}}
\newcommand{\UU}{{\mathcal U}}
\newcommand{\oboundary}{{\mathbb S}^1_\infty}
\newcommand{\Q}{{\mathbb Q}}
\newcommand{\DD}{{\mathcal D}}
\newcommand{\rot}{{\rm rot}}
\newcommand{\Cl}{{\rm Cl}}
\newcommand{\Index}{{\rm Index}}
\newcommand{\Int}{{\rm Int}}\newcommand{\Fix}{{\rm Fix}}
\newcommand{\Fr}{{\rm Fr}}
\newcommand{\ZZ}{\Z[2^{-1}]}
\newcommand{\II}{{\mathcal I}}
\newcommand{\JJ}{{\mathcal J}}
\newcommand{\g}{{\rm gen}}
\newcommand{\KK}{{\mathcal K}}
\newcommand{\OO}{{\mathcal O}}
\newcommand{\MM}{{\mathcal M}}
\newcommand{\NN}{{\mathcal N}}
\newcommand{\oM}{{\rho_M}}
\newcommand{\oG}{{H}}
\newcommand{\oH}{{\overline {\mathcal H}}}
\newcommand{\og}{{\overline g}}
\newcommand{\pre}{{\rm pre}}
\newcommand{\w}{{\underline{w}}}

\newcommand{\G}{{\mathbb G}}
\newcommand{\M}{\mathbb M}

\date{\today }

\maketitle

\begin{abstract}
We give a bijection between the isolated circular orders
of the group $G=PSL(2,\Z)\approx (\Z/2\Z)*(\Z/3\Z)$ and the equivalence classes of
Markov systems associated to $G$.
As applications, we present examples of isolated circular orders of the group
$G$.
\end{abstract}

\section{Introduction}

Throughout this paper, $G$ always stands for the group $(\Z/2\Z)*(\Z/3\Z)$.
Our purpose is
to give a bijection between the isolated circular orders
of the group $G$ and the symbolic dynamics associated with $G$, called
Markov systems.
We also apply it for
the construction of examples of isolated circular orders of 
$G$.
The paper is divided into four parts. In Part I, after preparing
necessary prerequisites, we define Markov systems of
 $G$ and state the main theorem.
Part II is devoted to the proof of one half of the main theorem,
and Part III of the other half.
In Part IV,  some examples of isolated circular orders are given.
In \cite{Ma}, 
 the space $LO(B_3)$ of the left orders (left invariant total orders) of the braid group
$B_3$ of three strings are shown to be homeomorphis to the space $CO(G)$ of the circular orders of
$G$. This provides examples of isolated left orders of $B_3$.
\bigskip

\begin{center}
 \large \sc Part I
\end{center}
After preparing necessary prerequisites in Sections 2--4,
 we  define Markov systems associated with $G$ and
state the main result  (Theorem \ref{T}) in Section 5.

\section{Circular orders}

In this section, we provide preliminary facts about circular orders.
Let $H$ be an arbitrary countable group.

\begin{definition}\label{co}
 A map $c:H^3\to\{0,1,-1\}$ is called a
{\em circular order of $H$} if it satisfies the following 
conditions

(1) $c(g_1,g_2,g_3)=0$ if and only if $g_i=g_j$ for some $i\ne j$.

(2) For any $g_1,g_2,g_3,g_4\in H$, we have
$$c(g_2,g_3,g_4)-c(g_1,g_3,g_4)+c(g_1,g_2,g_4)-c(g_1,g_2,g_3)=0.$$

(3) For any $g_1,g_2,g_3, g_4\in H$, we have
$$
c(g_4g_1,g_4g_2,g_4g_3)=c(g_1,g_2,g_3).$$
\end{definition}

\begin{definition}
 Given a finite 
set $F$ of $H$, a {\em configuration of $F$} is
an equivalence class of injections $\iota:F\to S^1$, where two
 injections $\iota$ and $\iota'$ is said to be equivalent if there is
an orientation preserving homeomorphism $h$ of $S^1$ such that $\iota'=h\iota$.
\end{definition}

Given a   circular order $c$ of $H$, the
configuration of the set $\{g_1,g_2,g_3\}$ of three points
is determined by the rule that $g_1,g_2,g_3$ is ordered
anticlockwise if $c(g_1,g_2, g_3)=1$, and clockwise if $c(g_1,g_2,g_3)=-1$.
By condition (2) of Definition \ref{co}, this is independent
of the enumeration of the set.
But (2) says more: an easy induction shows the following.
\begin{proposition}\label{should}
 A  circular order of $H$ determines the configuration of any
 finite subset $F$ of $\oG$. 
\end{proposition}

Denote by $CO(H)$ the set of all the   circular orders.
It is a closed subset of the space of maps from $\oG^3$ to $\{0,\pm1\}$,
and therefore
equipped with a totally disconnected compact metrizable topology.
 A circular order is said to be {\em isolated} if it is an isolated point
of $CO(H)$.
If $c\in CO(H)$ is isolated, then there is a finite subset
$S$ of $H$ such
that any circular order  which gives the same configuration of $S$
as $c$ is $c$, and {\em vice versa}.
Such a set $S$ is called a {\em determining set} of $c$.

For an automorphism $\sigma$ and $c\in CO(\oG)$, we define 
$\sigma_*c\in CO(\oG)$ by 
$$
(\sigma_*c)(g_1,g_2,g_3)=c(\sigma^{-1}g_1,\sigma^{-1}g_2,\sigma^{-1}g_3).
$$
The order $\sigma_*c$ is called an {\em automorphic image} of $c$.
An automorphic image of an isolated circular order is isolated.
We also say that $c$ and $\sigma_*c$ belong to the same {\em
automorphism class}.

Given $c\in CO(H)$, we define an action of $\oG$ on $S^1$
as follows.
Fix an enumeration of $H$: $H=\{g_i\mid i\in\N\}$ such that $g_1=e$
and a base point $x_0\in S^1$.
Define an embedding $\iota:H\to S^1$ inductively as follows. 
First, set $\iota(g_1)=x_0$ and $\iota(g_2)=x_0+1/2$. If $\iota$ is
defined on $\{g_1,\cdots,g_n\}$, then there is a connected component of 
$S^1\setminus\{\iota(g_1),\ldots,\iota(g_n)\}$ where the point $g_{n+1}$
should be embedded, by virtue of Proposition \ref{should}. Define
$\iota(g_{n+1})$ to be the midpoint of that interval.
The left translation of $\oG$ yields an action of $\oG$ on $\iota(\oG)$
which is shown to extend to a continuous action on $\Cl(\iota(\oG))$. Extend it
further to an action on $S^1$ by setting that the action on the
gaps\footnote{A gap of a closed subst $K$ of $S^1$ means a connected
component of $S^1\setminus K$.} be
linear.
The action so obtained is called the dynamical realization of $c$ based at $x_0$.

\begin{definition}
 An action $\phi$ of the group $\oG$ on $S^1$ is called {\em tight} at
$x_0\in S^1$ if it satisfies the following two conditions.

$\bullet$ It is free at $x_0$, i.e, the stabilizer of $x_0$ is trivial.

$\bullet$ If $J$ is a gap of the orbit
closure $\Cl(\phi(\oG)x_0)$, then $\partial J$
 is contained in the orbit $\phi(\oG)x_0$.
\end{definition}

We have three lemmas whose proofs are easy and omitted.
(Lemma \ref{tight} is a consequence of the midpoint construction.
In fact, it is used in the construction of the dynamical realization,
where we extend the $H$-action from the orbit of $x_0$ to the orbit closure.)

\begin{lemma}\label{tight}
 The dynamical realization is {\em tight} at the base point
$x_0$.
\end{lemma}

\begin{lemma}\label{conj}
 Two dynamical realizations obtained via different enumerations of $\oG$
 are mutually conjugate by an orientation and base point preserving
 homeomorphism of $S^1$.
\end{lemma}

\begin{lemma}\label{DR}
 An action $\phi$ of $\oG$ on $S^1$ which is tight at $x_0$ 
is topologically conjugate to the dynamical
realization of some circular order $c$ based at $x_0$
by an orientation and base point preserving homeomorphism.
\end{lemma}

 Henceforth any action $\phi$ as in the last lemma is
refered to as {\em a dynamical realization} of $c$.

\section{Preliminaries on $G$}

Here we study properties necessary for us of the group
$$
G=\langle\alpha,\beta\mid\alpha^2=\beta^3=e\rangle.
$$
As is well known, $G$ is isomorphic to $PSL(2,\Z)$, by an isomorphism $\phi$
such that 
$$
\phi(\alpha)=\left[\begin{array}{cc}0&-1\\1&0\end{array}\right],\ \ 
\phi(\beta)=\left[\begin{array}{cc}1&1\\ -1&0\end{array}\right].$$
See Figure 1 for the action of $PSL(2,\Z)$ on the Poincar\'e upper
half plane $\mathbb H$.
\begin{figure}[h]
{\unitlength 0.1in%
\begin{picture}( 20.8700, 20.9900)( 15.0000,-22.3000)%
%
\special{pn 4}%
\special{sh 1}%
\special{ar 2635 1269 16 16 0  6.28318530717959E+0000}%
\special{sh 1}%
\special{ar 2635 1269 16 16 0  6.28318530717959E+0000}%
%
\special{pn 8}%
\special{ar 2640 1269 947 947  5.6465485  5.6333485}%
%
\special{pn 8}%
\special{pa 2640 328}%
\special{pa 2640 2220}%
\special{fp}%
%
\special{pn 4}%
\special{sh 1}%
\special{ar 2145 1269 16 16 0  6.28318530717959E+0000}%
\special{sh 1}%
\special{ar 2145 1269 16 16 0  6.28318530717959E+0000}%
%
\special{pn 8}%
\special{pa 2640 323}%
\special{pa 2613 416}%
\special{pa 2603 447}%
\special{pa 2594 478}%
\special{pa 2574 540}%
\special{pa 2564 570}%
\special{pa 2554 601}%
\special{pa 2521 691}%
\special{pa 2509 721}%
\special{pa 2497 750}%
\special{pa 2484 779}%
\special{pa 2470 808}%
\special{pa 2456 836}%
\special{pa 2442 865}%
\special{pa 2427 892}%
\special{pa 2411 920}%
\special{pa 2377 974}%
\special{pa 2341 1026}%
\special{pa 2284 1104}%
\special{pa 2204 1204}%
\special{pa 2162 1254}%
\special{pa 2150 1269}%
\special{fp}%
%
\special{pn 8}%
\special{pa 2140 1259}%
\special{pa 2163 1282}%
\special{pa 2187 1305}%
\special{pa 2233 1351}%
\special{pa 2255 1374}%
\special{pa 2278 1398}%
\special{pa 2300 1421}%
\special{pa 2342 1469}%
\special{pa 2363 1494}%
\special{pa 2401 1544}%
\special{pa 2420 1570}%
\special{pa 2437 1596}%
\special{pa 2469 1650}%
\special{pa 2497 1706}%
\special{pa 2510 1735}%
\special{pa 2522 1764}%
\special{pa 2533 1793}%
\special{pa 2544 1823}%
\special{pa 2554 1854}%
\special{pa 2563 1884}%
\special{pa 2581 1946}%
\special{pa 2588 1978}%
\special{pa 2596 2009}%
\special{pa 2624 2137}%
\special{pa 2630 2170}%
\special{pa 2636 2202}%
\special{pa 2640 2220}%
\special{fp}%
%
\special{pn 8}%
\special{pa 2150 1269}%
\special{pa 1689 1269}%
\special{fp}%
%
\special{pn 4}%
\special{ar 2635 1269 117 117  4.2123938  1.2138844}%
%
\special{pn 4}%
\special{pa 2586 1171}%
\special{pa 2542 1191}%
\special{fp}%
\special{sh 1}%
\special{pa 2542 1191}%
\special{pa 2611 1182}%
\special{pa 2591 1169}%
\special{pa 2594 1145}%
\special{pa 2542 1191}%
\special{fp}%
%
\special{pn 4}%
\special{ar 2150 1264 165 165  3.4633432  5.0466612}%
%
\special{pn 4}%
\special{pa 1993 1200}%
\special{pa 1983 1230}%
\special{fp}%
\special{sh 1}%
\special{pa 1983 1230}%
\special{pa 2023 1173}%
\special{pa 2000 1179}%
\special{pa 1985 1160}%
\special{pa 1983 1230}%
\special{fp}%
\put(25.9200,-2.6100){\makebox(0,0)[lb]{$\infty$}}%
\put(15.0000,-13.1800){\makebox(0,0)[lb]{$-1$}}%
\put(20.8100,-13.7400){\makebox(0,0)[lb]{$\omega$}}%
\put(25.1500,-13.3200){\makebox(0,0)[lb]{$i$}}%
\put(28.6000,-13.3000){\makebox(0,0)[lb]{$\phi(\alpha)$}}%
\put(18.7000,-10.8000){\makebox(0,0)[lb]{$\phi(\beta)$}}%
\put(25.9000,-23.6000){\makebox(0,0)[lb]{$0$}}%
\put(23.9000,-12.5000){\makebox(0,0)[lb]{$P$}}%
\end{picture}}%

\caption{The open disk bounded by the circle is
the Poincar\'e upper half plane $\mathbb H$.
The element $\phi(\alpha)$ is the $1/2$-rotation
around $i$, and $\phi(\beta)$ the $1/3$-rotation around
$\omega=(-1+\sqrt{-3})/2$. 
The region $P$ bounded by the ideal triangle $\triangle 0\,\infty\,\omega$
is a fundamental domain of $PSL(2,\Z)$.
 }
\end{figure}
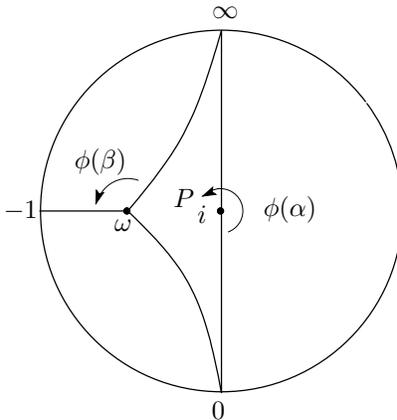

\begin{proposition}\label{p1}
(1) Any element of $G\setminus\{e\}$ is of order 2, 3 or infinite. 
Any element of order 2 is conjugate to $\alpha$, and any element of order
 3 is conjugate either to $\beta$ or to $\beta^{-1}$.

(2) Any torsion free subgroup of $G$ is isomorphic to a free group,
 either finitely generated or not.

(3) The commutator subgroup is a free group freely generated by 
$\alpha\beta\alpha\beta^{-1}$ and $\alpha\beta^{-1}\alpha\beta$.

(4) Any automorphism of $G$ is the conjugation by an element of 
$PGL(2,\Z)$ when we identify $G$ with $PSL(2,\Z)$ by $\phi$.
In other words, the outer automorphism group of $G$ is isomorphic to $\Z/2\Z$,
generated by the involutin $\sigma_0$ defined by
 $\sigma_0(\alpha)=\alpha$ and $\sigma_0(\beta)=\beta^{-1}$.
\end{proposition}

\bd (1) and (2) are well known, and can be obtained easily by
considering the action on $\mathbb H$.
(3) can be shown, for example, by an induction on word length of
elements in $[G,G]$, using the fact that the total exponents of $\alpha$
(resp.\ $\beta$) in
the word is even (resp.\ a multiple of 3).
 To show (4), notice that any automorphism $\sigma$ of $G$ sends $\beta$ to
an element which is conjugate either to $\beta$ or to $\beta^{-1}$.
We only need to verify that 
 $\sigma$ is an inner automorphism in the former case. By composing with
 an
appropriate inner
automorphism if necessary, we may assume that $\sigma(\alpha)=\alpha$. Let
$x$ be a fixed point of $\phi(\sigma(\beta))$ in $\mathbb H$. 
Since $\sigma(\beta)$ is a conjugate $\gamma\beta\gamma^{-1}$,
 $x$ is a translate of
$\omega$ in Figure 1 by an element $\phi(\gamma)\in PSL(2,\Z)$. All such points $x$,
except $\omega$ and $\phi(\alpha)(\omega)$, satisfy
$d(x,i)>d(\omega,i)$ where $d$ is the 
Poincar\'e distance, in which case $\phi(\alpha)$ and
$\phi(\sigma(\beta))$ generate a covolume infinite Fuchsian group. Since
 $\sigma$ is an automorphism, we have either
$x=\omega$ or $\phi(\alpha)\omega$. Accordingly,  $\sigma(\beta)=\beta$ or
$\sigma(\beta)=\alpha\beta\alpha$, In either case, $\sigma$ is an inner
automorphism.
\qed

\bigskip
{\sc Question}.
 Does the equality $(\sigma_0)_*(c)=-c$ hold? Or at least for isolated
 order
$c$?
This is true for all the examples constructed in Sections 13 and 14.

\section{Isolated circular orders of $G$}

Let  $c$  be an
isolated circular order of $G$, and $\rho$ a dynamical realization of $c$
based at $x_0\in S^1$.

\begin{proposition}
 \label{p31}
(1) There is a unique minimal set $\MM$ of the action $\rho$ 
homeomorphic to a Cantor set. Moreover $x_0\not\in\MM$.

(2) Let $I$ be a gap of $\MM$ which contains $x_0$. Then
all the gaps of $\MM$ is a translate of $I$ by the action $\rho$.
The stabilizer\footnote{$G_I=\{g\in G\mid\rho(g)I=I\}$} $G_I$
of $I$ is infinite cyclic. 

\end{proposition}

\bd It is known by \cite{MR} Corollary 1.3 that a dynamical
 realization
$\rho$ of an
isolated circular order cannot be minimal 
(for any countable group). Assume for contradiction that 
$\rho$ admits a finite minimal set of cardinality $n$. 
Considering the action of $G$ on the minimal set, one obtains a
 surjective homomorphism $\xi:G\to\Z/n\Z$.
Since $G/[G,G]$ is isomorphic to $\Z/6\Z$, $n$ is either 1, 2, 3
or 6. The circular order $c$ of $G$ induces a left
order $\lambda$ on $\rm Ker(\xi)$ by considering the action of
 $\rho\vert_{{\rm Ker}(\xi)}$ on the gap containing $x_0$. In particular, 
$\rm  Ker(\xi)$ must be
torsion free. This eliminates the case $n=1,2,3$. 
To eliminate the case $n=6$, we claim
that the induced left order $\lambda$ on ${\rm Ker}(\xi)$ must be
isolated.
Any left order $\lambda'$ on $\rm Ker(\xi)$ together with the cyclic
order on $\Z/6\Z$ defines a cyclic order $c'$ of $G$ by the lexicographic
construction. Moreover the map $\lambda'\mapsto c'$ is injective and
continuous. 
Since $c$ is isolated, $\lambda$
must be isolated, showing
the claim.
But if $n=6$, $\rm Ker(\xi)=[G,G]$ is a free group on 2 generators,
 and does not admit
an isolated left order \cite{McC}.
This shows that the minimal set cannot be a finite set.
It must be a Cantor set. The uniqueness of
a Cantor minimal set $\MM$ is easy and well known.  

If $x_0\in\MM$,
the two boundary points of a gap of $\MM=\Cl(\rho(G)x_0)$ cannot be
from one orbit of the action, contradicting the tightness of the action
$\rho$ at $x_0$. This shows $x_0\not\in\MM$.

The same argument shows that there are no gaps of $\MM$ which is not a
translate of $I$. Such a gap might be a gap of $\Cl(\rho(G)x_0)$.

Finally let us show that $G_I$ is infinite cyclic.
First, $G_I$ is nontrivial. Assume it is trivial and consider an
interval delimited by $x_0$ and a point $y\in\MM\subset\Cl(\rho(G)x_0)$. 
Then $y$ cannot be a translate of $x_0$, again
contradicting the tightness. If $G_I$ is not infinite cyclic,
$G_I$ is a free group on more than one generators, which does not admit
an isolated left order (\cite{McC}). As before, the lexicographic construction
gives a contradiction. \qed

\bigskip
The subgroup $G_I$ in the last proposition is called the {\em linear
subgroup} of the isolated circular order of $c$, and is denoted by
$L_c$. For an automorphism $\sigma$ of $G$, we have
$L_{\sigma_*c}=\sigma(L_c)$. 
For a fixed isolated circular order $c$, consider the minimal word
length of a generator of $\sigma(L_c)$ where $\sigma$ ranges over
all the automorphims of $G$.
It is an even number, say $2k$.

\begin{definition}\label{deg}
 The number $k$ is called the {\em degree} of the isolated circular
 order $c$, and is denoted by ${\rm deg}(c)$.
\end{definition}

Clearly the degree is an automorphism class function. It is an odd
number, as is remarked in  Section 14.

\section{Markov system}

\begin{definition}
 A {\em Markov system} ${\mathbb M}=(a,b,[a],[b],[b^{-1}])$ consists of an
 orientation preserving involution $a$, a period three homeomorphism $b$
of $S^1$ and subsets $[a],[b],[b^{-1}]$ of $S^1$ which satisfy
the conditions (A)--(E) to be listed below. 
\end{definition}

\smallskip\em
(A) The sets
 $[a]$, $[b]$ and $[b^{-1}]$ are disjoint,
 each consisting of
 $k$ closed intervals for some $k\in\N$.
The number $k$ is called the {\em multiplicity} of the system $\mathbb M$. \rm

\rm
\smallskip
A connected component of $[a]$ (resp.\ $[b]$ and $[b^{-1}]$) is called an
 {\em $a$-interval} (resp.\ $b$- and $b^{-1}$-{\em interval}).
Denote $X=[a]\cup[b\cup[b^{-1}]$.

\smallskip \em
(B) Two $a$-intervals are not adjacent in $X$.
 Likewise for $b$- and
 $b^{-1}$-intervals.

\smallskip\rm

A gap of $X$ between an $a$-interval and a $b^{\pm1}$-interval is
called a {\em principal gap}. A gap between a $b$-interval and a
$b^{-1}$-interval is called a {\em complementary gap}.
A maximal interval which consists of $b$-intervals and complementary
gaps is  called a {\em$b$-block}.
Any $b$-interval is contained in a unique $b$-block.
 The $a$-intervals and the $b$-blocks are alternating in $S^1$ and
there are just $k$ $b$-blocks.
The union of $b$-blocks is denoted by $[[b]]$.
Let us continue the conditions, which are reminiscent of the action
of $\alpha$ and $\beta$ on (the first letter of) the words representing elements
of $G$.


\smallskip\em
(C) $a[a]=[[b]]$
\smallskip\rm

\noindent
This implies $a[[b]]=[a]$ and hence $a[b^{\pm1}]\subset[a]$.

\smallskip\em
(D) $b[a]=[b]$, $b[b]=[b^{-1}]$.

\smallskip\rm\noindent
This implies $b[b^{-1}]=[a]$ since $b^3=id$, and also
$b^{-1}[a]=[b^{-1}]$, $b^{-1}[b^{-1}]=[b]$ and $b^{-1}[b]=[a]$.

\smallskip\rm

A principal gap $J$ is always mapped to a principal gap by $a$, and exactly
one of $b$ and $b^{-1}$ maps $J$ to a principal gap, the other to a
complementary gap. Therefore the principal gaps consist of several cycles.
Our last condition is the following.

\smallskip\em
(E) The principal gaps are formed of one cycle.
\rm

\begin{definition}
 Two Markov systems
${\mathbb M}=(a,b,[a],[b],[b^{-1}])$ and ${\mathbb M}'=(a',b',[a'],[b'],[(b')^{-1}])$ are said
to be {\em equivalent} if there is an orientation preserving
homeomorphism of $S^1$ which sends
 $[a]$, $[b]$ and
 $[b^{-1}]$
to $[a']$, $[b']$ and $[(b')^{-1}]$, respectively.
\end{definition}

Our main result is the following.

\begin{Theorem}
 \label{T}
There is a bijection between the isolated
 circular orders of $G$ and the equivalence classes of the
Markov systems.
\end{Theorem}

We would like to emphasize that this is not a true classification theorem
of isolated circular orders: it is too difficult to classify 
Markov systems. What we can do so far is to present examples of Markov
systems, as in Part IV.

\bigskip
\begin{center}
 \sc\Large Part II
\end{center}
We construct  Markov partitions associated to isolated circular orders,
thereby showing one half of Theorem \ref{T}. 
The result of Part II is summerized as Theorem \ref{At}. 

\section{Further properties of isolated circular orders of $G$}

This section is devoted to the preparation of 
basic facts needed for the proof of the
following theorem.
Let $\phi$ be a dynamical realization of an arbitrary isolated circular
 order of $G$, based at $x_0$, and $\MM$ the minimal set of the action $\phi$.
\begin{theorem}
 \label{At}
 There is a Markov partition $(a,b,[a],[b],[b^{-1}])$
such that $a=\phi(\alpha)$, $b=\phi(\beta)$.
\end{theorem}

 We denote $\G=\phi(G)$.  The group $\G$ is isomorphic
to $G$ (since $\phi$ is free at $x_0$) and is generated by $a$ and $b$. 

\begin{definition}
 A word of letters $a,b^{\pm1}$ is called {\em admissible} if
it is reduced and contains no consecutive same letters.
\end{definition}

Any element $g\in \G\setminus\{id\}$ is expressed uniquely by an
admissible word denoted by $W(g)$,
whose length  by $\Abs{g}$.
The first letter of $W(g)$ is called the {\em prefix} of $g$ and is denoted by $\pre(g)$.
We put $\pre(id)=\emptyset$ for completeness. If $W(g)=t_1t_2\cdots t_r$
and $W(g')=t_it_{i+1}\cdots t_r$ for some $2\leq i\leq r$, 
$g'$ is called a {\em larva} of $g$. 
The group $\G$ acts freely at the base point $x_0$ and all the above terminologies
about $\G$ are carried over to the orbit $\G x_0$.

\smallskip\em
If $x=gx_0$ and $x'=g'x_0$ for some $g,g'\in\G\setminus\{id\}$, 
the word $W(x)$ is to be the word $W(g)$,
the length $\Abs{x}$ is to be $\Abs{g}$, 
the prefix $\pre(x)$ is $\pre(g)$, and $x'$
is said to be a larva of $x$ if $g'$ is a larva of $g$. 

\smallskip\rm

Choose a Riemannian metric on $S^1$ in such a way that the involution $a$
is an isometry. The distance of two points $x,y\in S^1$ is denoted by
$\abs{x-y}$. The length of an interval $J$ is denoted by $\abs{J}$.
As a corollary of Proposition \ref{p31}, we get the following.

\begin{corollary}
 \label{Ac1}
There is $\epsilon_1>0$ such that if a closed interval $J$ satisfies
 $\abs{J}<\epsilon_1$ and $x_0\in J$, then $J\cap\G x_0=\{x_0\}$.
\end{corollary}

\bd Choose $\epsilon_1$ to be smaller than the distance of $x_0$ to
the neighbouring points in $\G x_0$. \qed

\begin{lemma}\label{Al31}
There are finitely many gaps $I_1,\ldots,I_r$ of the minimal set $\MM$
with the following properties.

(1) For any gap $J\neq I_i$, the prefixes of all the points
of $\G x_0\cap J$ are the same.

(2) For $J=I_i$, there is an enumeration of points of $\G x_0\cap J$:
$$\G x_0\cap J=\{\ldots,x_{-2},x_{-1},x_0,x_1,x_2\ldots\}$$
in the anti-clockwise order of $S^1$ such
that $\pre(x_n)$ ($n<0$) are the same, and 
$\pre(x_n)$ ($n>0$) are the same.
\end{lemma}

\bd Denote by $I$ the gap of $\MM$ containing $x_0$ as before and
by $G_I$ the stabilizer of $I$ by the action $\phi$.
Given an arbitrary gap $J$ of $\MM$, choose $g\in \G$ such that
$J=gI$ and that $\Abs{g}$ is the smallest among such $g$. 
Let $h$ be a generator of $\phi(G_I)$. Then the word $W(g)$
cannot have $W(h)$ or $W(h^{-1})$ as a larva.
This implies that if $\Abs{g}\geq\Abs{h}$,
 $W(g)$ is not completely
cancelled out in the word $W(gh^n)$.
In particular, $\pre(gh^n)$
is the same for any $n\in\Z$, showing (1). In the remaining case
$\Abs{g}<\Abs{h}$,  $\pre(gh^n)$ ($n>0$) are the same, as
well as $\pre(gh^{-n})$ ($n>0$), finishing the proof of (2). \qed

\begin{corollary}
 \label{Ac31}
There is $\epsilon_2>0$ such that if  a closed interval $J$
satisfies $\abs{J}<\epsilon_2$ and if there are points 
$x_1,x_2\in J\cap\G x_0$ such that $\pre(x_1)\neq\pre(x_2)$, then
$\Int(J)\cap\MM\neq\emptyset$.
\end{corollary}

\bd Choose $\epsilon_2>0$ smaller than the distance of $x_0$ 
and $x_{\pm1}$ for each of the intervals $I_1,\cdots,I_r$ in Lemma \ref{Al31}. \qed

\bigskip
For $t\in\{a,b^{\pm1}\}$, let $\G_t$ be the set of those elements $g\in \G\setminus\{id\}$ such
that the last letter of $W(g)$ is $t$. 

\begin{lemma}
 \label{Al35}
We have an inclusion
$\MM\subset\Cl(\G_t x_0)\setminus\G_t x_0$ for each $t\in\{a,b^{\pm1}\}$.

\end{lemma}

\bd  It
suffices to show that the closed set $X=\Cl(\G_t x_0)\setminus \G_t x_0$
 is invariant by $\G$, since $\MM$ is the unique minimal set.
Given $x\in X$, there is a sequence $\{g_n\}$ in $\G_t$ such that
$g_nx_0\to x$ and $\Abs{g_n}\to\infty$.
For any $f\in \G$, we have $fg_n\in \G_t$
if $n$ is sufficiently large,
showing that $fx=\lim_{n\to \infty} fg_nx_0\in X$. \qed

\bigskip
\begin{definition}
 An action $\psi:G\to{\rm Homeo}_+(S^1)$ is said to be {\em
 $\epsilon$-close} to the dynamical realization $\phi$ if 
$\Abs{\psi(\alpha)-\phi(\alpha)}_{0}<\epsilon$ and
$\Abs{\psi(\beta^{\pm1})-\phi(\beta^{\pm1})}_{0}<\epsilon$,
where $\psi(g)-\phi(g)$ is a map from the abelian group $S^1$ to $S^1$,
and $\Abs{\cdot}_0$ denotes the supremun norm.
\end{definition}

\begin{lemma}
 \label{Al32}
There is $\epsilon_3>0$ with the following property: if 
$\psi:G\to{\rm  Homeo}_+ (S^1)$ is $\epsilon_3$-close to $\phi$,
then $\psi$ is free at $x_0$ and the circular order of $G$
determined by the orbit $\psi(G)x_0\subset S^1$ is the same as $c$.
\end{lemma}

\bd Let $S\subset G$ be a finite determining set of $c$.
One can choose $\epsilon_3>0$ so that if $\psi$ is $\epsilon_3$-near to
$\phi$, then the circular order of $\psi(S)x_0$ is the same as
$\phi(S)x_0$ and additionally $\psi(\alpha)$ and $\psi(\beta)$ is not the identity.
Assume the isotropy group $H$ of $\psi$ at $x_0$ is nontrivial.
Then $H$ is torsion free, since any torsion element $\gamma$ of $G$ is conjugate
to $\alpha$ or $\beta^{\pm1}$ and $\psi(\gamma)$ is fixed point free
by the additional condition of $\epsilon_3$.
Therefore $H$ is a free group (Proposition 3.1) and admits  a left order $\lambda$ and its
reciproxal $-\lambda$. On the other hand, the quotient $G/H$ admits a
left $G$-invariant circular order $c'$ determined by
the orbit $\psi(G)x_0$. Now $\pm\lambda$ and $c'$ determines two distinct
circular orders by the lexicographic constructions. This is contrary to
the definition of the deterimining set $S$. We have shown that
$\psi$ is free at $x_0$. The rest of the assertion follows again by the
definition of the determining set. \qed

\section{Continuity of $W(x)$ }

In this section, we show that the assignment 
$\G x_0\ni x\mapsto W(x)$ is continuous in some weak sense.
The argument here follows closely the proof of \cite{MR}, Proposition
4.7, but we need an elaboration since the group $G$ is not a free
group treated in \cite{MR}. 
   
Two letters from the set $\{a,b^{\pm1}\}$ is said to be {\em congruent}
if either they are the same or one is $b$ and the other is $b^{-1}$.
Choose $\epsilon$ so that
$0<\epsilon<\min\{\epsilon_1,\epsilon_2,\epsilon_3, 1/2\}$,
where $\epsilon_1,\epsilon_2,\epsilon_3$ are the constants
defined in the last section.

\begin{proposition}
 \label{Ap41}
If $x,y\in\G x_0$ satisfy
$\abs{x-y}<\epsilon$, then $\pre(x)$ and $\pre(y)$ are congruent.
\end{proposition}

This section is devoted to the proof of the above proposition by 
contradiction. We assume the following conditions throughout this section,
and we shall deduce a contradiction at the end of the section.

\smallskip\em
 There is a closed interval $J$
with the following properties. 

\smallskip
($\sharp$1) $\abs{J}<\epsilon$.

($\sharp$2) $\partial J=\{x_1,x_2\}\subset \G x_0$ and
$\pre(x_1)$ and $\pre(x_2)$
are not congruent.

\smallskip\rm\noindent
By Corollay \ref{Ac1}, ($\sharp$1) and ($\sharp$2) imply the following ($\sharp$3).

\smallskip\em
($\sharp$3) $x_0\not\in J$.

\smallskip\rm\noindent
An interval which satisfies ($\sharp$1)--($\sharp$3) is called
a {\em $\sharp$-interval}.

\begin{lemma}
 \label{Al42}
There is a $\sharp$-interval $J$ 
such that no larva of a point in $\partial J$ is contained in $J$.
\end{lemma}

\newcommand{\sint}{{$\sharp$-interval\ }}
\bd Let $J=[x,y]$ be a \sint
 such that $N(J)=\Abs{x}+\Abs{y}$ is the smallest.
Then $x$ has no larva congruent to $x$ 
and contained in $J$. Likewise for $y$.
If there are no larvae of $x$ and $y$ contained in $J$ at all, we are
done.
So assume one of them, say $x$, has
a larva (necessarily 
not congruent to $x$) contained in $J$. Choose the larva
$z$ of $x$ in $J$ which is the nearest to $x$. Then the interval $J'=[x,z]$
is a $\sharp$-interval.
Moreover there are no larvae of $x$ other than $z$ contained in $J'$.

Now since $z$ is a larva of $x$, $x$ and $z$ belong to the
same $\G_t x_0$ for some  $t\in\{a,b^{\pm}\}$.
Choose $s\in\{a,b^{\pm1}\}\setminus\{t\}$.
By Corollary \ref{Ac31} and Lemma \ref{Al35}, 
$\G_sx_0\cap J'$ is nonempty since $J'$ is a $\sharp$-interval. Choose a point  $u$ of minimal length
from $\G_sx_0\cap J'$.
If $\pre(u)$ is congruent to $\pre(x)$,
choose the interval $[u,z]$: if not, $[x,u]$. \qed

\bigskip
Finally let us show that Lemma \ref{Al42} leads to a contradiction.
To fix the idea, let $J$ in Lemma \ref{Al42} be such that
$J=[x,y]$ and $\pre(x)=a$ and
$\pre(y)=b^{\pm1}$. Choose an open interval $U\supset J$ such that 
$\abs{U}<\epsilon$ and that
all the larvae of $x$ and $y$, as well as $x_0$, are not contained in
$U$. Let $h$ be a homeomorphism of $S^1$ supported on $U$ such that 
$h$ sends $x$ to the opposite side of $y$ in $U$, and define an action
$\psi$ of $G$ by setting $\psi(\alpha)=h\phi(\alpha)h^{-1}$ and
$\psi(\beta)=\phi(\beta)$. 
Notice that $\psi(\alpha)=\phi(\alpha)$ except on
$U\sqcup \phi(\alpha)U$. 
($U$ and $\phi(\alpha)U$ is disjoint since
$\epsilon<1/2$ and $\phi(\alpha)$ is an isometry.) 
By Lemma \ref{Al32},
 the cyclic order of $G$ obtained by the orbit
$\psi(G)x_0$ must be the same as $c$, i.e, the one obtained from $\G
x_0$.
However inductions on the length of $f$ and $g$ show that
$$(\psi(f)x_0,\psi(g)x_0,x_0)=(hx,y,x_0)\ \mbox{ and }\
(\phi(f)x_0,\phi(g)x_0,x_0)=(x,y,x_0).$$
A contradiction.
This completes the proof of Proposition \ref{Ap41}.

\section{Continuity of $ W(x)$ -continued}

\newcommand{\W}{{\widehat W}}

Let ${\mathcal W}_\infty$ be the set of infinite words of letters
$a,b^{\pm1}$ and $\emptyset$ with the following properties:

(1) there is no consecutive appearence of $a$, $b$ and $b^{-1}$,

(2) $b^{\mp1}$ does not follow $b^{\pm1}$, and

(3) all the letters after $\emptyset$ are $\emptyset$.

By certain abuse,
 we denote  ${ W}(g)\in{\mathcal W}_\infty$ to be a finite word $W(g)$
followed by a sequence of $\emptyset$. Thus for example
${W}(id)=\emptyset\,\emptyset\,\emptyset\cdots$, and
${W}(ab)=ab\,\emptyset\,\emptyset\cdots$.
For $n>0$, the initial subword of ${W}(g)$ of length $n$
is  denoted  ${W}_n(g)$,
We also define ${W}(x)$ and ${W}_n(x)$ for a point
$x=gx_0\in\G x_0$ by ${W}(x)={W}(g)$ and
${W}_n(x)={W}_n(g)$.
As a consequence of Proposition \ref{Ap41}, we get:
\begin{proposition}
 \label{Ap51}
For any $n>0$, there is $\epsilon(n)>0$ such that if two points
$x,y\in\G x_0$ satisfy $\abs{x-y}<\epsilon(n)$, then
${W}_n(x)={W}_n(y)$.
\end{proposition}

\bd For any $\eta>0$, define $\delta(\eta)\in(0,\eta)$ so that if
$\abs{x-y}<\delta(\eta)$, then 
$\abs{b^{\pm1}x-b^{\pm1}y}<\eta$.
(Recall that $a$ is assumed to 
be an isometry.)
Define
$\epsilon(n)=\overbrace{\delta(\delta(\cdots\delta))}^{n}(\epsilon )$.
An induction on $n$ shows the proposition. \qed

\section{Construction of Markov system}

 Given $x\in\MM$,
choose a point $x_i\in\G x_0$ from the $\epsilon(i)/2$-neighbourhood of $x$.
Then ${W}_i(x_i)$ is independent of the choice of $x_i$.
Moreover ${W}_i(x_i)={W}_i(x_{j})$
if $i<j$.
Thus the seqence $\{{W}_i(x_i)\}$ stabilizes.

\begin{definition}
 For any $x\in \MM$, define a word ${W}(x)\in {\mathcal W}_\infty$
as  the limit of ${W}_i(x_i)$.
Also define $W_n(x)$ to be the
initial subword of length $n$ of $W(x)$.
\end{definition}

Notice that $W(x)$ is a word of letters $a$ and $b^{\pm1}$:
$\emptyset$ never shows up.
\begin{definition}
 For an admissible word $w$ of length $n$ of letters $a$, $b^{\pm1}$ ,
we define the subset $[w]$ of $S^1$ to
 be the union of the points $x\in\MM$ such that ${W}_n(x)=w$ and
the gaps $(x,y)$ of $\MM$ such that ${W}_n(x)={W}_n(y)=w$.
\end{definition}

\begin{lemma}
 \label{Al51}
(1) For any finite admissible word $w$, $[w]$ is a finite union of closed
 disjoint intervals.

(2) If $v\ne w$ are redued words of the same length, then
 $[v]\cap[w]=\emptyset$.

(3) We have
 $a[ab^{\pm1}]=[b^{\pm1}]$.

(4) We have $b[a]=[b]$,
$b[b]=[b^{-1}]$, and
 $b[b^{-1}]=[a]$.

(5) The cardinalities of the components of $[a]$, $[b]$ and $[b^{-1}]$
are the same.
\end{lemma}

\bd (1) is a consequence of Proposition \ref{Ap51}. 
(2)--(4) are clear from the definitions.
(5) follows from (4). \qed

\bigskip
{\sc Proof of Theorem \ref{At}.}
Define $\M=(a,b,[a],[b,[b^{-1}]])$.
It is obvious that $\M$ satisfies the conditions (A)--(E).
\qed


\bigskip
\begin{center}
 \sc\Large Part III
\end{center}
We define properties (*) and (**), and show that any Markov systems are
equivalent to one with (*) and (**).
Next, to Markov partitions with (*) and (**), we assign
isolated circular orders,
thereby showing the other half of Theorem \ref{T}. The result is summerized
as Theorem \ref{t1}.

\section{Fundamental properties of Markov systems and modifications}

Let $\M=(a,b,[a],[b],[b^{-1}])$ be an arbitrary Markov system.
Recall that the map $a$ sends a principal gap $J$ to a principal
gap, and either one of $b$ or $b^{-1}$ sends $J$ to a principal
gap, the other to a complementary gap.
By condition (E),
the principal gaps $I_i$, $I_i'$ and complementary gaps 
$I_i''$, $i\in \Z/k\Z$, are dynamically related as in (\ref{e1}) below, where $b_i$ is
either $b$ or $b^{-1}$. The gaps $I_i'''$ are gaps between components
of $[ab]$ and $[ab^{-1}]$ to be defined later.

\begin{equation}
 \label{e1}
\begin{array}{cccccccccc}
I_1&&\stackrel{b_1}{-\!\!\!\longrightarrow}&&I_1'\stackrel{a}{\to}I_2&&\stackrel{b_2}{-\!\!\!\longrightarrow}&&I_2'&
 \stackrel{a}{\to}\cdots\stackrel{b_{k}}{\to}\ \ \ I_{k}'\stackrel{a}{\to}I_{k+1}=I_1\\

&\stackrel{b_1}{\nwarrow}&&\stackrel{b_1}{\swarrow}&&
\stackrel{b_2}{\nwarrow}&&\stackrel{b_2}{\swarrow}&&\\
&&I''_1&&&&I''_2&&&\\
&&a\downarrow&&&&a\downarrow&&&\\
&&I_1'''&&&&I_2'''&&&
\end{array}
\end{equation}

In this section, we study fundamental properties of Markov systems.
Besides, we show that any Markov system can be modified 
in its equivalence class to
another one with good properties.
First of all, consider the following property concerning diagram (\ref{e1}).

\smallskip
\em
(*) The map\footnote{In $f_1=ab_k\cdots ab_1$,
$b_1$ is the first map.} $f_1=ab_k\cdots ab_1$ which leaves $I_1$ invariant
  admits no fixed point in the open interval
$I_1$, and for any $z\in I_1$,
$\displaystyle\lim_{n\to\infty}f^n_1(z)\to x$,
where $x$ is a point in $\partial J\cap[a]$.
\rm\smallskip

Our first modification result is the following.

\begin{lemma}
 \label{l68}
Any Markov system is equivalent to a Markov system with
(*).
\end{lemma}

\bd
In the admissible word $f_1=ab_{k}\cdots ab_1$, the last map $a$ (first
in the word) is a
transposition of $I'_{k}$ and $I_1$. If we
change 
the map $a$ by the conjugation by a homeomorphism supported in
$I_1$ and
leave $b$ unchanged, the new maps still satisfy all
the requirement for Markov systems.
 Since the modification is free, one can get (*). \qed

\bigskip

We introduce basic terminologies and notations.

\begin{definition}
 For subsets $P$ and $Q$ consisting of finite disjoint closed intervals
 of $S^1$,
the inclusion $P\subset Q$ is called {\em precise} if any boundary point
of $Q$ is contained in $P$.
\end{definition}
The inclusion $[b]\cup[b^{-1}]\subset[[b]]$ is precise.
The composite of precise inclusions is precise.

\begin{definition}
For the Markov system $\M=(a,b,[a],[b],[b^{-1}])$,
denote by $\G=\G(\M)$ the subgroup of ${\rm Homeo}_+(S^1)$ generated
by $a$ and $b$.
\end{definition}

By {\em word}, we always mean word of letters $a,b^{\pm1}$.
Any map of $\G\setminus\{id\}$ is expressed uniquely as an admissible
word, as can be shown by (4) and (7) of the next lemma.

\begin{definition}
  For a  admissible
word $w=vt$ where $t$ is the last letter of $w$, 
we define $[w]=v[t]$.
\end{definition}

\begin{lemma}
 \label{l61}
(1) For an admissible word $wv$, we have $w[v]=[wv]$.

(2) If $wv$, $v^{-1}u$ and $wu$ are admissible, then 
$wv[v^{-1}u]=[wu]$.

(3) If $wa$ is admissible, $[wa]=[w]$.

(4) If $wv$ is admissible, $[wv]\subset[w]$.

(5) For an admissible word $w$, $[w]$ consists of $k$ disjoint closed
 intervals.

(6) The inclusion $[wb]\cup[wb^{-1}]\subset[w]$ is precise for an
admissible word $w$ which ends at $a$.

(7) If $w$ and $w'$ are distinct admissible words of the same
 length, $[w]$ and $[w']$ are disjoint.
\end{lemma}

\bd (1) is obtained by an easy induction on the length of $v$. For (2),
$wv[v^{-1}u]=wvv^{-1}[u]=w[u]=[wu]$. For (3), if $w=vb^{\pm1}$,
$[wa]=vb^{\pm1}[a]=v[b^{\pm1}]=[vb^{\pm1}]=[w]$. For (4), notice that if
$w=ua$ and $v=b^{\pm1}$, $[wv]=ua[b^{\pm1}]\subset u[a]=[w]$. 
Together with (3), this imples $[wv]\subset[w]$ if $\Abs{v}=1$. 
The general case follows by an induction on $\Abs{v}$. For (5), if
$w=vt$ where $t$ is the last letter of $v$,  $[vt]=v[t]$, $[t]$
consists of disjoint $k$ intervals and $v$ is a homeomorphism.
For (6), since the inclusion
$[b]\cup[b^{-1}]\subset[[b]]$ is precise, if we write $w=va$, the inclusion
$$[wb]\cup[wb^{-1}]=w([b]\cup[b^{-1}])\subset w[[b]]=va[[b]]=v[a]=[w]$$
is precise. 
(7) follows from an induction of word length.
\qed

\bigskip

 For an infinite admissible word $\w=t_1t_2\cdots$, define 
$[\w]=\cap_{i}[t_1\cdots t_i]$. 
It is a nonempty set consisting
 of $k$ closed intervals, some possibly 
 degenerate to points.
We have $v[\w]=[v\w]$ for any finite admissible word $v$ and
any infinite admissible word $\w$.

\begin{definition}
Define $X_\infty=\cap_{f\in \G}\, f^{-1}X$ for $X=[a]\cup[b]\cup[b^{-1}]$.
\end{definition}

The set $X_\infty$ is closed and $\G$-invariant.

\begin{lemma}
 \label{l62}
We have $X_\infty=\cup_\w\,[\w]$, where $\w$ runs over all the infinite
 admissible 
words.
\end{lemma}

\bd The inclusion $\supset$ is easy: 
for any  $f\in \G$, we have $f[\w]=[f\w]\subset X$,
showing $[\w]\subset f^{-1}X$.
Let us show $\subset$.
For any $x\in X_\infty$, 
  define $t_1\in\{a,b^{\pm1}\}$ by the condition $x\in[t_1]$, then
$t_2$ by $t_1^{-1}(x)\in[t_2]$,  $t_3$ by
$t_2^{-1}t_1^{-1}(x)\in[t_3]$ $\ldots$. The word $\w=t_1t_2\cdots$ we obtained
is admissible and $x\in[\w]$. \qed

\bigskip

\em For any $x\in X_\infty$, define
${\hat W}(x)=\w$ if $x\in[\w]$.
\rm
 The inclusion
$X_\infty\subset X$ is ``precise'' in the sense that any boundary point
of $X$ is contained in $X_\infty$. This stems from Lemma \ref{l61} (6)
and Lemma \ref{l62}.
Therefore a gap of $X$ is a gap of $X_\infty$.

\begin{lemma}
 \label{l63}
The gaps of $X_\infty$ form one orbit of the $\G$-action.
\end{lemma}

\bd Any gap of $X_\infty$ other than principal or complementary gaps is
a gap contained in $[w]$
between components of $[wb]$ and  $[wb^{-1}]$ for
some admissible word $w$ ending at $a$, and hence is the
image by $w$ of a complementary gap. \qed

\begin{lemma}
 \label{Cantor}
The set $X_\infty$ admits no isolated component.
\end{lemma}

\bd The proof is by contradiction.
Let $C$ be an isolated component contained in
$[\w]$. Let $\w=t_1t_2t_3\cdots$ and for each $m\in\N$, let
$C_m$ be the component of
$[t_1\cdots t_m]$ containing $C$. We claim that the decreasing sequence
$C_1\supset C_2\supset\cdots\downarrow C$ stabilizes, i.e, there is $m_0$ such that
$C_{m_0+i}=C_{m_0}$ for any $i\in\N$. Assume not. For any
neighbourhood $U$ of $C$, some $C_m$ is contained in $U$.
But if the sequence does not stabilize, there is a component of
$[\w]$ distinct from $C$ contained in $C_m$ and hence in $U$. 
This shows that $C$ is not isolated. The contradiction shows the
claim.
Now the interval $t_{m_0-1}^{-1}\cdots t_{1}^{-1}C$ is at the same time a
component of $X$ and of $X_\infty$, 

It is no loss of generality to assume 
that $C$ itself is a component of both $X$ and
$X_\infty$. 
Let $C(1)=C$, and $C(i)=t_{i-1}^{-1}\cdots t_{1}^{-1}C$ for any $i>1$.
 Then $C(i)$ is also a component of both $X$ and $X_\infty$.
Since $X$ has only finitely many components, the sequence $\{C(i)\}$ is
eventually periodic.

Without loosing generality,
one may assume that $C(1)$ is in a periodic cycle: $C(n+1)=C(1)$ and $n$
such smallest. Notice that the intertwining arrows of the cycle are
$t_i^{-1}:C(i)\to C(i+1)$ and that $C(i)$ is a $t_i$-interval. 
In particular, 

\smallskip
(1) if $C(i)$ is an $a$-interval,  $t_i^{-1}=a$.

\smallskip\noindent
Let $J(i)$ be the gap of $X$ right to $C(i)$. 
It is also a gap of $X_\infty$, either principal or complementary.
They form a cycle with the same arrows $t_i^{-1}:J(i)\to J(i+1)$ as 
the arrows $t_i^{-1}:C(i)\to C(i+1)$ of the cycle $\{C(i)\}$.
The cycle, consisting of principal and complementary gaps,
is contained in the first and second lines of diagram (\ref{e1}).
But since  there is no consecutive $b$ or $b^{-1}$ in 
the arrows of the cycle, it is the cycle in the first line or its reciprocal.
 In particular, $n=2k$ and all the $J(i)$'s are principal.

 Now any $a$-interval appears in the cycle $\{C(i)\}$,
since its right gap is  principal and any principal gap appears in
the cycle $\{J(i)\}$.
By (1), if $C(i)$ is an $a$-interval, it is mapped by $a$ to a
$b^{\pm1}$- interval $C(i+1)$. That is, {\em any $a$-interval must be mapped by
$a$ to a $b^{\pm1}$-interval}.
But  there are $k$ $a$-intervals 
and $2k$ $b^{\pm1}$-intervals, and some $a$-interval must be mapped by
$a$ to a nontrivial $b$-block containing more than one
$b^{\pm1}$-intervals. A contradiction. \qed

\bigskip

We consider the following property and the second modification.

\smallskip
\em
(**) ${\rm Int}X_\infty=\emptyset$.
\rm

\begin{lemma}
 \label{l66} Any Markov system is equivalent to a system
with properties (*) and (**).
\end{lemma}

\bd This is done by an anti-Denjoy modification:
one collapses each component $[\w]$ of $X_\infty$
to a point, and define new maps $a$ and $b$ of the collapsed $S^1$. 
The modification does not collapse $[a]$- or $[b^{\pm1}]$-intervals 
to points, thanks to Lemma \ref{Cantor}, and does not spoil property (*). \qed

\section{Circular order determined by Markov system}

In this section, we shall show that any Markov system 
with (*) and (**) is a dynamical realization of some 
circular order, i.e, that the associated action is tight at some point.
 This is done by showing that the set $X_\infty$ is a minimal set
with good properties.

\begin{definition}
 For an admissible word $v$ of even length, $(v)$ denotes the infinite
admissible word which
repeats $v$.
\end{definition}

\begin{lemma}
 \label{l64}
The boundary point $x\in\partial I_1\cap[a]$ satisfies 
${\hat W}(x)=(ab_{k}\cdots ab_1)$ and the other boundary point $y$ satisfies
${\hat W}(y)=(b_1^{-1}a\cdots b_{k}^{-1}a)$.
\end{lemma}

\bd Since $x$ is fixed by $f_1$, we have an equality of
infinite words
$$
ab_{k}\cdots ab_1{\hat W}(x)={\hat W}(x),$$
where the LHS is before reducing. If there is no cancelation in the LHS,
 we get ${\hat W}(x)=(ab_{k}\cdots ab_1)$.
If $ab_{k}\cdots ab_1$ is completely canceled out, we get 
${\hat W}(x)=(b_1^{-1}a\cdots b_{k}^{-1}a)$.
Otherwise there is an intermediate cancelation,
and the LHS begins at
$a$ and the RHS begins at $b_1^{-1}$.
We thus obtained either ${\hat W}(x)=(ab_k\cdots ab_1)$ or
${\hat W}(x)=(b_1^{-1}a\cdots b_k^{-1}a)$.
Like statement holds for ${\hat W}(y)$.
But since $x\in[a]$, ${\hat W}(x)=(ab_{k}\cdots ab_1)$, 
and since $I_1$ is a principal gap,
 ${\hat W}(y)=(b_1^{-1}a\cdots b_{k}^{-1}a)$. \qed

\begin{lemma}
 \label{l65} The stabilizer $\G_{I_1}$ of $I_1$ is infinite cyclic
generated by $f_1$. 
\end{lemma}

\bd Assume $h\neq id$ stabilizes $I_1$.
One may assume that $h$ 
admits neither $f_1$ nor
$f_1^{-1}$ 
as an initial subword, 
by replacing $h$ with
a shorter word if necessary.
Since $hx=x$ for $x$ in the previous lemma,
 we have an equality of infinite word
\begin{equation}
 \label{e2}
h(ab_{k}\cdots ab_1)=(ab_{k}\cdots ab_1).
\end{equation}
If there is an intermediate cancelation in the LHS, then $h$ begins and
ends at $a$. But then since $hy=y$, we have another equality
\begin{equation}
\label{e3}
h(b_1^{-1}a\cdots b_{k}^{-1}a)=(b_1^{-1}a\cdots b_{k}^{-1}a),
\end{equation}
where the LHS begins at $a$ and the RHS at $b_1^{-1}$, leading to a
contradiction.

Consider the case where there is no cancelation at all in the LHS of
(\ref{e2}).
Then 
 either $h$ is $f_1$ or its
initial subword of length $2\ell$, $\ell<k$.
In the latter case, the word $(ab_k\cdots ab_1)$ is periodic of period
$2\ell$ and $2k$, hence of period $2(k,\ell)$.
This shows that $h=ab_\ell\cdots ab_1$.
But then in the diagram
(\ref{e1}), the principal gap $I_{\ell+1}$ must be equal to $I_1$.
A contradiction.
In the remaining case where $w$ is cancelled out in the LHS of (\ref{e2}),
there is no cancelation in the LHS of (\ref{e3}), and one can show
that $h=f_1^{-1}$ by a like argument. \qed

\begin{lemma}
 \label{l67}
For any Markov system with (**),
 $X_\infty$ is a minimal set of $\G$.
\end{lemma}

\bd 
First of all, notice that $X_\infty$
is a Cantor set by Lemma \ref{Cantor} and (**). This, together with Lemma \ref{l63},
shows that the orbit of a boundary
point of a gap is dense in $X_\infty$. 
Assume there is a minimal set $Y$ properly contained in $X_\infty$.
Then any boundary point of a gap of $X_\infty$ cannot be contained in
$Y$. Therefore any gap $K$ of $Y$ contains infinitely many gaps $J_i$ of
$X_\infty$. Since each  $J_i$ belongs to the orbit of $I_1$, it
is left invariant by a map represented in an admissible word as
  $h_i=v_if_iv_i^{-1}$, where $f_i$ is a cyclic permutation of $f_1$.
Since all the $h_i$ leaves a boundary point $z$ of $K$ invariant, 
we get $v_if_iv_i^{-1}{\hat W}(z)={\hat W}(z)$.
By the same argument as the proof of Lemma \ref{l64},
${\hat W}(z)=v_i(f_i)$ or
${\hat W}(z)=v_i(f_i^{-1})$ for any $i$.
This contradicts the uniqueness  of ${\hat W}(z)$:
 the admissibility of the word $v_if_iv_i^{-1}$ implies
that $v_i$ contains neither $f_i$ nor $f_i^{-1}$ as the  terminal
subword.
\qed

\begin{definition}
 Given a Markov system
$\M=(a,b,[a],[b],[[b]])$,
 define a homomorphism $\phi_\M:G\to{\rm Homeo}_+(S^1)$ 
by $\phi_\M(\alpha)=a$ and $\phi_\M(\beta)=b$. 
\end{definition}

We have $\phi_\M(G)=\G$.

\begin{lemma}
 \label{l69}
For any Markov system $\M$ with (*) and (**),  $\phi_\M$ is
a dynamical realization of a circular order, denoted by $c_\M$,
 based at some point
 $x_0\in I_1$.
\end{lemma}

\bd 
By Lemma \ref{DR}, we only need to show that $\phi_\M$ is tight at $x_0$.
Clearly $\phi_\M$ is free at $x_0$ by (*). To show the other condition,
notice that $X_\infty$ is a Cantor minimal set and that any gap of
$X_\infty$ is a translate of $I_1$.
Now choose an arbitrary gap $J$ of $\Cl(\phi_\M(G) x_0)$ for $x_0\in I_1$. If one  end point of
$J$ belongs to the orbit $\phi_\M(G) x_0$, so does the other end point.
 On the other hand, there is no
gap of $\Cl(\phi_\M(G) x_0)$ whose end points belong to the minimal
set $X_\infty$. \qed

\section{Isolation}

The following theorem is the goal of Part II.

\begin{theorem}
 \label{t1}
For any Markov system $\M$ with (*) and (**),
the circular order $c_\M$ given by Lemma \ref{l69} is isolated.
Moreover ${\rm deg}(c_\M)$ is equal to the multiplicity of $\M$.
\end{theorem}

\bigskip
The statement about the degree (Definition \ref{deg}) is obvious. 
The rest of this section is devoted to the proof of the isolation
 of $c_\M$, by a pingpong argument.
Let 
$$Y=[ab]\cup[ab^{-1}]\cup[b]\cup[b^{-1}],\ \ 
h_1=abab^{-1},\  h_2=ab^{-1}ab,$$ 
$$\Omega(h_1)=[ab],\ \ \ \Omega(h_2)=[ab^{-1}], \ \ \
\Omega(h_1^{-1})=[b],\mbox{ and } \Omega(h_2^{-1})=[b^{-1}].$$ 
Notice that the ``address'' of $\Omega(h_i^{\pm1})$ is just the first two
letters of $h_i^{\pm1}$, since $[b^{\pm1}]=[b^{\pm1}a]$.

\begin{lemma}
 \label{l610}
For $i=1,2$, we have a precise inclusion
$$
h_i(Y-\Omega(h_i^{-1}))\subset\Omega(h_i)\ \ \mbox{ and }\ \ 
h_i^{-1}(Y-\Omega(h_i))\subset\Omega(h_i^{-1}).
$$
\end{lemma}

\bd For example, the first inclusion for $i=1$ follows from a sequence of precise
inclusions;
$$
abab^{-1}([a]\cup[b^{-1}])=aba([b^{-1}]\cup[b])
=ab([ab^{-1}]\cup[ab])\subset
ab([a])=[ab].$$
The other inclusions are similar. \qed

\bigskip
However the pingpoing property on $Y$ is not sufficient for
the proof of the isolation of $c_\M$.
\begin{lemma}
 \label{l611}
There are open neighbourhoods $N_i^{\pm}$ of $\Omega(h_i^{\pm1})$
 such that

(1) the closures $\Cl(N_i^{\pm})$ are mutually disjoint, and

(2)
$h_i(S^1-N_i^{-})\subset N_i^+$, or equivalently,
$h_i^{-1}(S^1-N_i^+)\subset N_i^{-}$, $i=1,2$.
\end{lemma}

\bd To fix the idea, we assume that the point $x$ in $\partial I_1\cap[a]$ is
the right end point of $I_1$.
 Thus the map $f_1=ab_k\cdots ab_1$ is an
increasing homeomorphism of $I_1$. 
One can choose two points in all the
principal gaps $I_i$ and $I_i'$ in (\ref{e1}) which satisfy the relations in Figure 2,
where the thin lines 
indicate the correspondences by the intertwining maps.
\begin{figure}[h]
{\unitlength 0.1in%
\begin{picture}( 45.9900, 10.9000)(  9.1000,-16.9000)%
%
\special{pn 20}%
\special{pa 1000 790}%
\special{pa 1000 1690}%
\special{fp}%
%
\special{pn 20}%
\special{pa 2422 817}%
\special{pa 2422 1690}%
\special{fp}%
%
\special{pn 20}%
\special{pa 2773 799}%
\special{pa 2773 1672}%
\special{fp}%
%
\special{pn 20}%
\special{pa 3340 799}%
\special{pa 3340 1672}%
\special{fp}%
%
\special{pn 20}%
\special{pa 3682 781}%
\special{pa 3682 1672}%
\special{fp}%
%
\special{pn 20}%
\special{pa 5140 799}%
\special{pa 5140 1663}%
\special{fp}%
%
\special{pn 20}%
\special{pa 5500 790}%
\special{pa 5500 1663}%
\special{fp}%
%
\special{pn 4}%
\special{sh 1}%
\special{ar 1000 1015 16 16 0  6.28318530717959E+0000}%
\special{sh 1}%
\special{ar 1000 1015 16 16 0  6.28318530717959E+0000}%
%
\special{pn 4}%
\special{sh 1}%
\special{ar 2422 1033 16 16 0  6.28318530717959E+0000}%
\special{sh 1}%
\special{ar 2773 1033 16 16 0  6.28318530717959E+0000}%
\special{sh 1}%
\special{ar 3340 1060 16 16 0  6.28318530717959E+0000}%
\special{sh 1}%
\special{ar 3682 1042 16 16 0  6.28318530717959E+0000}%
\special{sh 1}%
\special{ar 5149 1069 16 16 0  6.28318530717959E+0000}%
\special{sh 1}%
\special{ar 5500 1069 16 16 0  6.28318530717959E+0000}%
\special{sh 1}%
\special{ar 5500 1078 16 16 0  6.28318530717959E+0000}%
%
\special{pn 4}%
\special{sh 1}%
\special{ar 1000 1447 16 16 0  6.28318530717959E+0000}%
\special{sh 1}%
\special{ar 2422 1474 16 16 0  6.28318530717959E+0000}%
\special{sh 1}%
\special{ar 2773 1465 16 16 0  6.28318530717959E+0000}%
\special{sh 1}%
\special{ar 3340 1474 16 16 0  6.28318530717959E+0000}%
\special{sh 1}%
\special{ar 3682 1483 16 16 0  6.28318530717959E+0000}%
\special{sh 1}%
\special{ar 5140 1474 16 16 0  6.28318530717959E+0000}%
\special{sh 1}%
\special{ar 5500 1465 16 16 0  6.28318530717959E+0000}%
\special{sh 1}%
\special{ar 5500 1465 16 16 0  6.28318530717959E+0000}%
%
\special{pn 4}%
\special{pa 991 1456}%
\special{pa 2422 1033}%
\special{fp}%
%
\special{pn 4}%
\special{pa 2773 1033}%
\special{pa 2422 1492}%
\special{fp}%
%
\special{pn 4}%
\special{pa 2782 1456}%
\special{pa 2935 1303}%
\special{fp}%
%
\special{pn 4}%
\special{pa 3340 1051}%
\special{pa 3214 1150}%
\special{fp}%
%
\special{pn 4}%
\special{pa 3673 1042}%
\special{pa 3358 1474}%
\special{fp}%
%
\special{pn 4}%
\special{pa 3682 1483}%
\special{pa 5140 1069}%
\special{fp}%
%
\special{pn 4}%
\special{pa 5509 1078}%
\special{pa 5140 1474}%
\special{fp}%
\put(9.1000,-7.4500){\makebox(0,0)[lb]{$I_1$}}%
\put(23.3200,-7.5400){\makebox(0,0)[lb]{$I_i'$}}%
\put(27.2800,-7.3600){\makebox(0,0)[lb]{$I_2$}}%
\put(32.3200,-7.3600){\makebox(0,0)[lb]{$I'_{k-1}$}}%
\put(37.0000,-7.3600){\makebox(0,0)[lb]{$I_k$}}%
\put(50.5000,-7.3600){\makebox(0,0)[lb]{$I'_k$}}%
\put(53.9000,-7.3000){\makebox(0,0)[lb]{$I_1$}}%
\end{picture}}%

\caption{The right side of the intervals are drawn upward.}
\end{figure}
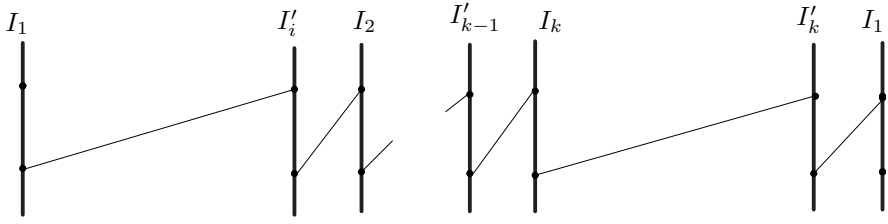
This is possible simply because the first return map
 is increasing.

 Then choose
three points in  $I_i''$ and $I_i'''$ as the images of the previous
 points. See Figure 3.
All these points are called
{\em distinguished points}.
\begin{figure}[h]
{\unitlength 0.1in%
\begin{picture}( 45.9900, 19.9500)(  9.1000,-25.9500)%
%
\special{pn 20}%
\special{pa 1000 790}%
\special{pa 1000 1690}%
\special{fp}%
%
\special{pn 20}%
\special{pa 2422 817}%
\special{pa 2422 1690}%
\special{fp}%
%
\special{pn 20}%
\special{pa 2773 799}%
\special{pa 2773 1672}%
\special{fp}%
%
\special{pn 20}%
\special{pa 3340 799}%
\special{pa 3340 1672}%
\special{fp}%
%
\special{pn 20}%
\special{pa 3682 781}%
\special{pa 3682 1672}%
\special{fp}%
%
\special{pn 20}%
\special{pa 5140 799}%
\special{pa 5140 1663}%
\special{fp}%
%
\special{pn 20}%
\special{pa 5500 790}%
\special{pa 5500 1663}%
\special{fp}%
%
\special{pn 20}%
\special{pa 2026 1312}%
\special{pa 2026 2266}%
\special{fp}%
%
\special{pn 20}%
\special{pa 4762 1330}%
\special{pa 4762 2230}%
\special{fp}%
%
\special{pn 4}%
\special{sh 1}%
\special{ar 1000 1015 16 16 0  6.28318530717959E+0000}%
\special{sh 1}%
\special{ar 1000 1015 16 16 0  6.28318530717959E+0000}%
%
\special{pn 4}%
\special{sh 1}%
\special{ar 2422 1033 16 16 0  6.28318530717959E+0000}%
\special{sh 1}%
\special{ar 2773 1033 16 16 0  6.28318530717959E+0000}%
\special{sh 1}%
\special{ar 3340 1060 16 16 0  6.28318530717959E+0000}%
\special{sh 1}%
\special{ar 3682 1042 16 16 0  6.28318530717959E+0000}%
\special{sh 1}%
\special{ar 5149 1069 16 16 0  6.28318530717959E+0000}%
\special{sh 1}%
\special{ar 5500 1069 16 16 0  6.28318530717959E+0000}%
\special{sh 1}%
\special{ar 5500 1078 16 16 0  6.28318530717959E+0000}%
%
\special{pn 4}%
\special{sh 1}%
\special{ar 1000 1447 16 16 0  6.28318530717959E+0000}%
\special{sh 1}%
\special{ar 2422 1474 16 16 0  6.28318530717959E+0000}%
\special{sh 1}%
\special{ar 2773 1465 16 16 0  6.28318530717959E+0000}%
\special{sh 1}%
\special{ar 3340 1474 16 16 0  6.28318530717959E+0000}%
\special{sh 1}%
\special{ar 3682 1483 16 16 0  6.28318530717959E+0000}%
\special{sh 1}%
\special{ar 5140 1474 16 16 0  6.28318530717959E+0000}%
\special{sh 1}%
\special{ar 5500 1465 16 16 0  6.28318530717959E+0000}%
\special{sh 1}%
\special{ar 5500 1465 16 16 0  6.28318530717959E+0000}%
%
\special{pn 4}%
\special{sh 1}%
\special{ar 2026 1519 16 16 0  6.28318530717959E+0000}%
\special{sh 1}%
\special{ar 2026 1789 16 16 0  6.28318530717959E+0000}%
\special{sh 1}%
\special{ar 2026 2059 16 16 0  6.28318530717959E+0000}%
\special{sh 1}%
\special{ar 2026 2059 16 16 0  6.28318530717959E+0000}%
%
\special{pn 4}%
\special{sh 1}%
\special{ar 4762 1519 16 16 0  6.28318530717959E+0000}%
\special{sh 1}%
\special{ar 4762 1807 16 16 0  6.28318530717959E+0000}%
\special{sh 1}%
\special{ar 4762 2095 16 16 0  6.28318530717959E+0000}%
\special{sh 1}%
\special{ar 4762 2086 16 16 0  6.28318530717959E+0000}%
%
\special{pn 4}%
\special{pa 1018 1042}%
\special{pa 2035 1537}%
\special{fp}%
%
\special{pn 4}%
\special{pa 1000 1447}%
\special{pa 2026 1789}%
\special{fp}%
%
\special{pn 4}%
\special{pa 991 1456}%
\special{pa 2422 1033}%
\special{fp}%
%
\special{pn 4}%
\special{pa 2431 1033}%
\special{pa 2026 1789}%
\special{fp}%
%
\special{pn 4}%
\special{pa 2422 1465}%
\special{pa 2026 2068}%
\special{fp}%
%
\special{pn 4}%
\special{pa 2773 1033}%
\special{pa 2422 1492}%
\special{fp}%
%
\special{pn 4}%
\special{pa 2782 1456}%
\special{pa 2935 1303}%
\special{fp}%
%
\special{pn 4}%
\special{pa 3340 1051}%
\special{pa 3214 1150}%
\special{fp}%
%
\special{pn 4}%
\special{pa 3673 1042}%
\special{pa 3358 1474}%
\special{fp}%
%
\special{pn 4}%
\special{pa 3682 1483}%
\special{pa 5140 1069}%
\special{fp}%
%
\special{pn 4}%
\special{pa 3682 1492}%
\special{pa 4771 1807}%
\special{fp}%
%
\special{pn 4}%
\special{pa 5140 1078}%
\special{pa 4780 1798}%
\special{fp}%
%
\special{pn 4}%
\special{pa 5509 1078}%
\special{pa 5140 1474}%
\special{fp}%
%
\special{pn 4}%
\special{pa 5140 1474}%
\special{pa 4762 2104}%
\special{fp}%
%
\special{pn 20}%
\special{pa 1342 1699}%
\special{pa 1342 2581}%
\special{fp}%
%
\special{pn 20}%
\special{pa 4042 1708}%
\special{pa 4042 2581}%
\special{fp}%
%
\special{pn 4}%
\special{pa 2008 1528}%
\special{pa 1342 1942}%
\special{fp}%
%
\special{pn 4}%
\special{pa 2026 1798}%
\special{pa 1342 2194}%
\special{fp}%
%
\special{pn 4}%
\special{pa 2026 2050}%
\special{pa 1351 2464}%
\special{fp}%
%
\special{pn 4}%
\special{pa 3682 1042}%
\special{pa 4762 1519}%
\special{fp}%
%
\special{pn 4}%
\special{pa 4762 1528}%
\special{pa 4042 1951}%
\special{fp}%
%
\special{pn 4}%
\special{pa 4762 1798}%
\special{pa 4042 2185}%
\special{fp}%
%
\special{pn 4}%
\special{pa 4762 2086}%
\special{pa 4033 2428}%
\special{fp}%
\put(9.1000,-7.4500){\makebox(0,0)[lb]{$I_1$}}%
\put(23.3200,-7.5400){\makebox(0,0)[lb]{$I_i'$}}%
\put(27.2800,-7.3600){\makebox(0,0)[lb]{$I_2$}}%
\put(32.3200,-7.3600){\makebox(0,0)[lb]{$I'_{k-1}$}}%
\put(37.0000,-7.3600){\makebox(0,0)[lb]{$I_k$}}%
\put(50.5000,-7.3600){\makebox(0,0)[lb]{$I'_k$}}%
\put(11.9800,-27.2500){\makebox(0,0)[lb]{$I_1'''$}}%
\put(19.3600,-24.2800){\makebox(0,0)[lb]{$I''_1$}}%
\put(39.6100,-27.0700){\makebox(0,0)[lb]{$I'''_k$}}%
\put(46.8100,-23.7400){\makebox(0,0)[lb]{$I''_k$}}%
\put(53.9000,-7.3000){\makebox(0,0)[lb]{$I_1$}}%
\end{picture}}%

\caption{}
\end{figure}
Define the neighbourhood $N_i^{\pm}$ of $\Omega(h_i^{\pm1})$ by
expanding each component until it reachs the first distinguised point.
Condition (1) of the lemma is satisfied. 

Let us check if (2) is satisfied. First of all,
we only consider the right gap $J$ of each component of $\Omega(h_i^{\mp1})$
and check if $J\setminus N_i^{\mp}$ is mapped by $h_i^{\pm1}$ into 
$h_i^{\pm1}(J)\cap N_i^{\pm}$. 
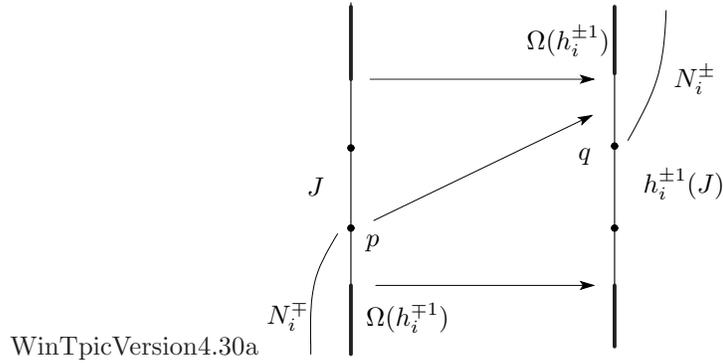
\begin{figure}[h]
WinTpicVersion4.30a
{\unitlength 0.1in%
\begin{picture}( 21.3000, 18.4000)( 11.6000,-22.8000)%
%
\special{pn 8}%
\special{pa 1600 440}%
\special{pa 1600 2250}%
\special{fp}%
%
\special{pn 8}%
\special{pa 2980 460}%
\special{pa 2980 2240}%
\special{fp}%
%
\special{pn 20}%
\special{pa 1600 460}%
\special{pa 1600 840}%
\special{fp}%
%
\special{pn 20}%
\special{pa 2980 460}%
\special{pa 2980 810}%
\special{fp}%
%
\special{pn 20}%
\special{pa 1600 1920}%
\special{pa 1600 2280}%
\special{fp}%
%
\special{pn 20}%
\special{pa 2980 1920}%
\special{pa 2980 2240}%
\special{fp}%
%
\special{pn 4}%
\special{pa 1700 840}%
\special{pa 2870 840}%
\special{fp}%
\special{sh 1}%
\special{pa 2870 840}%
\special{pa 2803 820}%
\special{pa 2817 840}%
\special{pa 2803 860}%
\special{pa 2870 840}%
\special{fp}%
%
\special{pn 4}%
\special{pa 1730 1920}%
\special{pa 2870 1920}%
\special{fp}%
\special{sh 1}%
\special{pa 2870 1920}%
\special{pa 2803 1900}%
\special{pa 2817 1920}%
\special{pa 2803 1940}%
\special{pa 2870 1920}%
\special{fp}%
%
\special{pn 4}%
\special{pa 1720 1580}%
\special{pa 2860 1030}%
\special{fp}%
\special{sh 1}%
\special{pa 2860 1030}%
\special{pa 2791 1041}%
\special{pa 2812 1053}%
\special{pa 2809 1077}%
\special{pa 2860 1030}%
\special{fp}%
%
\special{pn 4}%
\special{pa 1530 1650}%
\special{pa 1513 1678}%
\special{pa 1481 1734}%
\special{pa 1466 1762}%
\special{pa 1438 1820}%
\special{pa 1427 1849}%
\special{pa 1417 1879}%
\special{pa 1409 1910}%
\special{pa 1403 1941}%
\special{pa 1398 1972}%
\special{pa 1394 2004}%
\special{pa 1390 2068}%
\special{pa 1389 2101}%
\special{pa 1389 2198}%
\special{pa 1390 2231}%
\special{pa 1390 2280}%
\special{fp}%
%
\special{pn 4}%
\special{pa 3050 1160}%
\special{pa 3066 1132}%
\special{pa 3082 1103}%
\special{pa 3097 1075}%
\special{pa 3127 1017}%
\special{pa 3141 989}%
\special{pa 3155 960}%
\special{pa 3168 930}%
\special{pa 3180 901}%
\special{pa 3191 871}%
\special{pa 3201 841}%
\special{pa 3210 811}%
\special{pa 3224 749}%
\special{pa 3229 718}%
\special{pa 3234 686}%
\special{pa 3238 654}%
\special{pa 3244 590}%
\special{pa 3248 526}%
\special{pa 3249 493}%
\special{pa 3250 480}%
\special{fp}%
%
\special{pn 4}%
\special{sh 1}%
\special{ar 2980 1190 16 16 0  6.28318530717959E+0000}%
\special{sh 1}%
\special{ar 2980 1190 16 16 0  6.28318530717959E+0000}%
%
\special{pn 4}%
\special{sh 1}%
\special{ar 2980 1620 16 16 0  6.28318530717959E+0000}%
\special{sh 1}%
\special{ar 2980 1620 16 16 0  6.28318530717959E+0000}%
%
\special{pn 4}%
\special{sh 1}%
\special{ar 1600 1200 16 16 0  6.28318530717959E+0000}%
\special{sh 1}%
\special{ar 1600 1200 16 16 0  6.28318530717959E+0000}%
%
\special{pn 4}%
\special{sh 1}%
\special{ar 1600 1620 16 16 0  6.28318530717959E+0000}%
\special{sh 1}%
\special{ar 1600 1620 16 16 0  6.28318530717959E+0000}%
\put(16.8000,-21.7000){\makebox(0,0)[lb]{$\Omega(h_i^{\mp1})$}}%
\put(11.6000,-21.6000){\makebox(0,0)[lb]{$N_i^{\mp}$}}%
\put(32.9000,-9.2000){\makebox(0,0)[lb]{$N_i^{\pm}$}}%
\put(13.7000,-14.5000){\makebox(0,0)[lb]{$J$}}%
\put(31.3000,-14.7000){\makebox(0,0)[lb]{$h_i^{\pm1}(J)$}}%
\put(25.2000,-7.2000){\makebox(0,0)[lb]{$\Omega(h_i^{\pm1})$}}%
\put(16.8000,-17.4000){\makebox(0,0)[lb]{$p$}}%
\put(27.9000,-12.9000){\makebox(0,0)[lb]{$q$}}%
\end{picture}}%

\caption{The arrows are $h_i^{\pm1}$. We shall check if $p$ is mapped above
 $q$.
If this is true, the inverse $h_i^{\mp1}$ maps $q$ below $p$.}
\end{figure}
We do not need to consider the left gaps, since the statement for them
follows automatically. We can see this by Figure 4.
Our strategy is the following. We choose any
gap $J\in\{I_i,I'_i,I''_i,I'''_i\}$ from diagram (\ref{e1}). 
We consider the left end point $\partial_-J$
 and calculate to which class ($[ab]$, $[ab^{-1}]$, $[b]$ or
$[b^{-1}]$) it belongs. It tells us
which one of $\Omega(h_i^{\pm1})$ the left neighbour of $J$ is,
and thus which one of the maps $h_i^{\pm1}$ 
we should check.
We explain it concretely with an example $J=I_1$.
Recall that the left end point $y$ of $I_1$ satisfies 
${\hat W}(y)=(b_1^{-1}a\cdots b_k^{-1}a)$ and it belongs to 
$[b_1^{-1}]=\Omega(b_1^{-1}ab_1a)$. Therefore the map we should check is
$ab_1^{-1}ab_1$. For any gap $J$, we calculate its left end point and
the map in concern this way. 
The actual proof is divided into four cases.

\smallskip\noindent
{\sc Case 1}. $J=I_i$. Since 
$${\hat W}(\partial_-I_i)=ab_{i-1}\cdots ab_1{\hat W}(y)
=ab_{i-1}\cdots ab_1(b_1^{-1}a\cdots b_k^{-1}a)=b_i^{-1}ab_{i+1}^{-1}\cdots,
$$
we have $\partial_-I_i\in[b_i^{-1}]=\Omega(b_i^{-1}ab_ia)$ and the map in
concern is $h=ab_i^{-1}ab_i$. 

According as $b_{i+1}=b_i^{-1}$ or 
$b_{i+1}=b_i$, $h$ is
either of the following composite.
$$
 \begin{array}{cc}
I_i\stackrel{b_i}{\to}I_i'\stackrel{a}{\to}I_{i+1}&\stackrel{b_{i+1}}
{\to} I_{i+1}'\stackrel{a}{\to}I_{i+2}\\
 &\stackrel{b_{i+1}^{-1}}{\searrow} I_{i+1}''\stackrel{a}{\to}I_{i+1}'''
\end{array}
$$

In any case, $h$ satisfies the required property. See Figure 5.
\begin{figure}[h]

{\unitlength 0.1in%
\begin{picture}( 42.6000, 20.0000)( 10.5000,-24.0000)%
%
\special{pn 4}%
\special{ar 1490 590 42 42  5.2528085  5.4977871}%
%
\special{pn 4}%
\special{pa 1210 420}%
\special{pa 1210 1430}%
\special{fp}%
%
\special{pn 4}%
\special{pa 1990 400}%
\special{pa 1990 1400}%
\special{fp}%
%
\special{pn 4}%
\special{pa 2790 420}%
\special{pa 2790 1420}%
\special{fp}%
%
\special{pn 4}%
\special{pa 4390 420}%
\special{pa 4390 1420}%
\special{fp}%
%
\special{pn 4}%
\special{pa 5210 430}%
\special{pa 5210 1380}%
\special{fp}%
%
\special{pn 4}%
\special{sh 1}%
\special{ar 1210 700 16 16 0  6.28318530717959E+0000}%
\special{sh 1}%
\special{ar 1210 700 16 16 0  6.28318530717959E+0000}%
%
\special{pn 4}%
\special{sh 1}%
\special{ar 1210 1180 16 16 0  6.28318530717959E+0000}%
\special{sh 1}%
\special{ar 1210 1170 16 16 0  6.28318530717959E+0000}%
%
\special{pn 4}%
\special{sh 1}%
\special{ar 1990 700 16 16 0  6.28318530717959E+0000}%
\special{sh 1}%
\special{ar 1990 700 16 16 0  6.28318530717959E+0000}%
%
\special{pn 4}%
\special{sh 1}%
\special{ar 1990 1180 16 16 0  6.28318530717959E+0000}%
\special{sh 1}%
\special{ar 1990 1180 16 16 0  6.28318530717959E+0000}%
%
\special{pn 4}%
\special{sh 1}%
\special{ar 2780 700 16 16 0  6.28318530717959E+0000}%
\special{sh 1}%
\special{ar 2780 700 16 16 0  6.28318530717959E+0000}%
%
\special{pn 4}%
\special{sh 1}%
\special{ar 2790 1200 16 16 0  6.28318530717959E+0000}%
\special{sh 1}%
\special{ar 2790 1200 16 16 0  6.28318530717959E+0000}%
%
\special{pn 4}%
\special{sh 1}%
\special{ar 4390 710 16 16 0  6.28318530717959E+0000}%
\special{sh 1}%
\special{ar 4390 710 16 16 0  6.28318530717959E+0000}%
%
\special{pn 4}%
\special{sh 1}%
\special{ar 5210 720 16 16 0  6.28318530717959E+0000}%
\special{sh 1}%
\special{ar 5200 720 16 16 0  6.28318530717959E+0000}%
%
\special{pn 4}%
\special{sh 1}%
\special{ar 5210 1210 16 16 0  6.28318530717959E+0000}%
\special{sh 1}%
\special{ar 5210 1220 16 16 0  6.28318530717959E+0000}%
%
\special{pn 4}%
\special{sh 1}%
\special{ar 4390 1210 16 16 0  6.28318530717959E+0000}%
\special{sh 1}%
\special{ar 4390 1210 16 16 0  6.28318530717959E+0000}%
%
\special{pn 8}%
\special{pa 1270 1150}%
\special{pa 1920 760}%
\special{fp}%
\special{sh 1}%
\special{pa 1920 760}%
\special{pa 1853 777}%
\special{pa 1874 787}%
\special{pa 1873 811}%
\special{pa 1920 760}%
\special{fp}%
%
\special{pn 8}%
\special{pa 2060 1170}%
\special{pa 2730 740}%
\special{fp}%
\special{sh 1}%
\special{pa 2730 740}%
\special{pa 2663 759}%
\special{pa 2685 769}%
\special{pa 2685 793}%
\special{pa 2730 740}%
\special{fp}%
%
\special{pn 8}%
\special{pa 2880 1180}%
\special{pa 4290 750}%
\special{fp}%
\special{sh 1}%
\special{pa 4290 750}%
\special{pa 4220 750}%
\special{pa 4239 766}%
\special{pa 4232 789}%
\special{pa 4290 750}%
\special{fp}%
%
\special{pn 8}%
\special{pa 4460 1180}%
\special{pa 5140 760}%
\special{fp}%
\special{sh 1}%
\special{pa 5140 760}%
\special{pa 5073 778}%
\special{pa 5095 788}%
\special{pa 5094 812}%
\special{pa 5140 760}%
\special{fp}%
%
\special{pn 8}%
\special{pa 3590 1240}%
\special{pa 3590 2180}%
\special{fp}%
%
\special{pn 8}%
\special{pa 4210 1580}%
\special{pa 4210 1580}%
\special{pa 4210 2400}%
\special{pa 4210 2400}%
\special{pa 4210 1580}%
\special{pa 4210 1580}%
\special{fp}%
%
\special{pn 4}%
\special{sh 1}%
\special{ar 3590 1460 16 16 0  6.28318530717959E+0000}%
\special{sh 1}%
\special{ar 3590 1460 16 16 0  6.28318530717959E+0000}%
%
\special{pn 4}%
\special{sh 1}%
\special{ar 3590 1740 16 16 0  6.28318530717959E+0000}%
\special{sh 1}%
\special{ar 3590 1740 16 16 0  6.28318530717959E+0000}%
%
\special{pn 4}%
\special{sh 1}%
\special{ar 3590 2020 16 16 0  6.28318530717959E+0000}%
\special{sh 1}%
\special{ar 3590 2020 16 16 0  6.28318530717959E+0000}%
%
\special{pn 4}%
\special{sh 1}%
\special{ar 4210 1730 16 16 0  6.28318530717959E+0000}%
\special{sh 1}%
\special{ar 4210 1730 16 16 0  6.28318530717959E+0000}%
%
\special{pn 4}%
\special{sh 1}%
\special{ar 4210 1960 16 16 0  6.28318530717959E+0000}%
\special{sh 1}%
\special{ar 4210 1960 16 16 0  6.28318530717959E+0000}%
%
\special{pn 4}%
\special{sh 1}%
\special{ar 4210 2190 16 16 0  6.28318530717959E+0000}%
\special{sh 1}%
\special{ar 4210 2190 16 16 0  6.28318530717959E+0000}%
%
\special{pn 4}%
\special{pa 2880 740}%
\special{pa 3510 1440}%
\special{fp}%
\special{sh 1}%
\special{pa 3510 1440}%
\special{pa 3480 1377}%
\special{pa 3474 1400}%
\special{pa 3451 1404}%
\special{pa 3510 1440}%
\special{fp}%
%
\special{pn 4}%
\special{pa 3680 1490}%
\special{pa 4120 1700}%
\special{fp}%
\special{sh 1}%
\special{pa 4120 1700}%
\special{pa 4068 1653}%
\special{pa 4072 1677}%
\special{pa 4051 1689}%
\special{pa 4120 1700}%
\special{fp}%
\put(10.5000,-12.4000){\makebox(0,0)[lb]{$p$}}%
\put(53.1000,-7.8000){\makebox(0,0)[lb]{$q$}}%
\put(43.3000,-17.7000){\makebox(0,0)[lb]{$q'$}}%
\end{picture}}%

\caption{$h$ maps $p$ above $q$ or $q'$.} 
\end{figure}
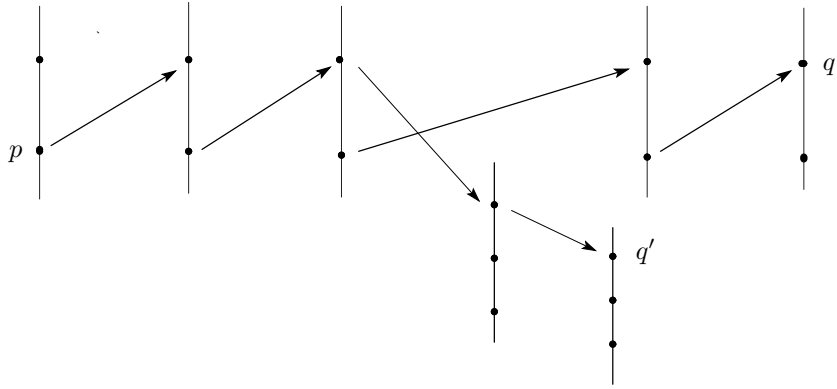

\smallskip\noindent
{\sc Case 2}. $J=I_i'$.
Since 
${\hat W}(\partial_-I_i')=b_i{\hat W}(\partial_-I_i)=ab_{i+1}^{-1}\cdots$,
we have $\partial_-I_i'\in\Omega(ab_{i+1}^{-1}ab_{i+1})$ and the map in
concern is $b_{i+1}^{-1}ab_{i+1}a$. This case is analogous to Case 1,
and is omitted. The point is that the map $b_{i+1}^{-1}ab_{i+1}a$ defined on $I_i'$ also
goes from left to right in (\ref{e1}).

\smallskip\noindent
{\sc Case 3}. $J=I_i''$.
Since 
${\hat W}(\partial_-I_i'')=b_i^{-1}{\hat W}(\partial_-I_i)=b_i\cdots$,
we have $\partial_-I_i''\in\Omega(b_iab_i^{-1}a)$ and the map in
concern is $ab_iab_i^{-1}$.

 It is either of the following composite.
$$
 \begin{array}{ccc}
 &  I_i'\stackrel{a}{\to}{I_{i+1}}&\stackrel{b_{i+1}}{\to}I'_{i+1}\stackrel{a}
{\to}I_{i+2}\\
I_i''\stackrel{b_i^{-1}}{\nearrow}&
 &\stackrel{b_{i+1}^{-1}}{\searrow}I_{i+1}''
\stackrel{a}{\to}I_{i+1}''' \\
\end{array}
$$
See Figure 6. In the bottom composite, the point $p$ is mapped
exactly to $q'$. In this case, we enlarge the neighbourhood of
$\Omega(b_iab_i^{-1}a)$ slightly so as to contain the point $p$.
Then its complement in $J$ is mapped above $q'$.
\begin{figure}[h]
{\unitlength 0.1in%
\begin{picture}( 43.4000, 20.8000)(  9.7000,-24.8000)%
%
\special{pn 4}%
\special{ar 1490 590 42 42  5.2528085  5.4977871}%
%
\special{pn 4}%
\special{pa 1990 400}%
\special{pa 1990 1400}%
\special{fp}%
%
\special{pn 4}%
\special{pa 2790 420}%
\special{pa 2790 1420}%
\special{fp}%
%
\special{pn 4}%
\special{pa 4390 420}%
\special{pa 4390 1420}%
\special{fp}%
%
\special{pn 4}%
\special{pa 5210 430}%
\special{pa 5210 1380}%
\special{fp}%
%
\special{pn 4}%
\special{sh 1}%
\special{ar 1990 700 16 16 0  6.28318530717959E+0000}%
\special{sh 1}%
\special{ar 1990 700 16 16 0  6.28318530717959E+0000}%
%
\special{pn 4}%
\special{sh 1}%
\special{ar 1990 1180 16 16 0  6.28318530717959E+0000}%
\special{sh 1}%
\special{ar 1990 1180 16 16 0  6.28318530717959E+0000}%
%
\special{pn 4}%
\special{sh 1}%
\special{ar 2780 700 16 16 0  6.28318530717959E+0000}%
\special{sh 1}%
\special{ar 2780 700 16 16 0  6.28318530717959E+0000}%
%
\special{pn 4}%
\special{sh 1}%
\special{ar 2790 1200 16 16 0  6.28318530717959E+0000}%
\special{sh 1}%
\special{ar 2790 1200 16 16 0  6.28318530717959E+0000}%
%
\special{pn 4}%
\special{sh 1}%
\special{ar 4390 710 16 16 0  6.28318530717959E+0000}%
\special{sh 1}%
\special{ar 4390 710 16 16 0  6.28318530717959E+0000}%
%
\special{pn 4}%
\special{sh 1}%
\special{ar 5210 720 16 16 0  6.28318530717959E+0000}%
\special{sh 1}%
\special{ar 5200 720 16 16 0  6.28318530717959E+0000}%
%
\special{pn 4}%
\special{sh 1}%
\special{ar 5210 1210 16 16 0  6.28318530717959E+0000}%
\special{sh 1}%
\special{ar 5210 1220 16 16 0  6.28318530717959E+0000}%
%
\special{pn 4}%
\special{sh 1}%
\special{ar 4390 1210 16 16 0  6.28318530717959E+0000}%
\special{sh 1}%
\special{ar 4390 1210 16 16 0  6.28318530717959E+0000}%
%
\special{pn 8}%
\special{pa 2060 1170}%
\special{pa 2730 740}%
\special{fp}%
\special{sh 1}%
\special{pa 2730 740}%
\special{pa 2663 759}%
\special{pa 2685 769}%
\special{pa 2685 793}%
\special{pa 2730 740}%
\special{fp}%
%
\special{pn 8}%
\special{pa 2880 1180}%
\special{pa 4290 750}%
\special{fp}%
\special{sh 1}%
\special{pa 4290 750}%
\special{pa 4220 750}%
\special{pa 4239 766}%
\special{pa 4232 789}%
\special{pa 4290 750}%
\special{fp}%
%
\special{pn 8}%
\special{pa 4460 1180}%
\special{pa 5140 760}%
\special{fp}%
\special{sh 1}%
\special{pa 5140 760}%
\special{pa 5073 778}%
\special{pa 5095 788}%
\special{pa 5094 812}%
\special{pa 5140 760}%
\special{fp}%
%
\special{pn 8}%
\special{pa 3590 1240}%
\special{pa 3590 2180}%
\special{fp}%
%
\special{pn 8}%
\special{pa 4210 1580}%
\special{pa 4210 1580}%
\special{pa 4210 2400}%
\special{pa 4210 2400}%
\special{pa 4210 1580}%
\special{pa 4210 1580}%
\special{fp}%
%
\special{pn 4}%
\special{sh 1}%
\special{ar 3590 1460 16 16 0  6.28318530717959E+0000}%
\special{sh 1}%
\special{ar 3590 1460 16 16 0  6.28318530717959E+0000}%
%
\special{pn 4}%
\special{sh 1}%
\special{ar 3590 1740 16 16 0  6.28318530717959E+0000}%
\special{sh 1}%
\special{ar 3590 1740 16 16 0  6.28318530717959E+0000}%
%
\special{pn 4}%
\special{sh 1}%
\special{ar 3590 2020 16 16 0  6.28318530717959E+0000}%
\special{sh 1}%
\special{ar 3590 2020 16 16 0  6.28318530717959E+0000}%
%
\special{pn 4}%
\special{sh 1}%
\special{ar 4210 1730 16 16 0  6.28318530717959E+0000}%
\special{sh 1}%
\special{ar 4210 1730 16 16 0  6.28318530717959E+0000}%
%
\special{pn 4}%
\special{sh 1}%
\special{ar 4210 1960 16 16 0  6.28318530717959E+0000}%
\special{sh 1}%
\special{ar 4210 1960 16 16 0  6.28318530717959E+0000}%
%
\special{pn 4}%
\special{sh 1}%
\special{ar 4210 2190 16 16 0  6.28318530717959E+0000}%
\special{sh 1}%
\special{ar 4210 2190 16 16 0  6.28318530717959E+0000}%
%
\special{pn 4}%
\special{pa 2880 740}%
\special{pa 3510 1440}%
\special{fp}%
\special{sh 1}%
\special{pa 3510 1440}%
\special{pa 3480 1377}%
\special{pa 3474 1400}%
\special{pa 3451 1404}%
\special{pa 3510 1440}%
\special{fp}%
%
\special{pn 4}%
\special{pa 3680 1490}%
\special{pa 4120 1700}%
\special{fp}%
\special{sh 1}%
\special{pa 4120 1700}%
\special{pa 4068 1653}%
\special{pa 4072 1677}%
\special{pa 4051 1689}%
\special{pa 4120 1700}%
\special{fp}%
\put(53.1000,-7.8000){\makebox(0,0)[lb]{$q$}}%
\put(43.3000,-17.7000){\makebox(0,0)[lb]{$q'$}}%
%
\special{pn 4}%
\special{pa 1210 1260}%
\special{pa 1210 2150}%
\special{fp}%
%
\special{pn 4}%
\special{sh 1}%
\special{ar 1210 1460 16 16 0  6.28318530717959E+0000}%
\special{sh 1}%
\special{ar 1210 1460 16 16 0  6.28318530717959E+0000}%
%
\special{pn 4}%
\special{sh 1}%
\special{ar 1210 1740 16 16 0  6.28318530717959E+0000}%
\special{sh 1}%
\special{ar 1210 1740 16 16 0  6.28318530717959E+0000}%
%
\special{pn 4}%
\special{sh 1}%
\special{ar 1210 1980 16 16 0  6.28318530717959E+0000}%
\special{sh 1}%
\special{ar 1210 1970 16 16 0  6.28318530717959E+0000}%
%
\special{pn 4}%
\special{pa 1300 1950}%
\special{pa 1900 1280}%
\special{fp}%
\special{sh 1}%
\special{pa 1900 1280}%
\special{pa 1841 1316}%
\special{pa 1864 1320}%
\special{pa 1870 1343}%
\special{pa 1900 1280}%
\special{fp}%
\put(9.7000,-20.5000){\makebox(0,0)[lb]{$p$}}%
\put(11.1000,-11.8000){\makebox(0,0)[lb]{$I_i''$}}%
\put(15.6000,-18.4000){\makebox(0,0)[lb]{$b_i^{-1}$}}%
%
\special{pn 20}%
\special{pa 1210 2150}%
\special{pa 1210 2480}%
\special{fp}%
\put(13.5000,-23.7000){\makebox(0,0)[lb]{$\Omega(b_iab_i^{-1}a)$}}%
\end{picture}}%

\caption{$p$ is mapped abouve $q$, but exactly to $q'$.}
\end{figure}
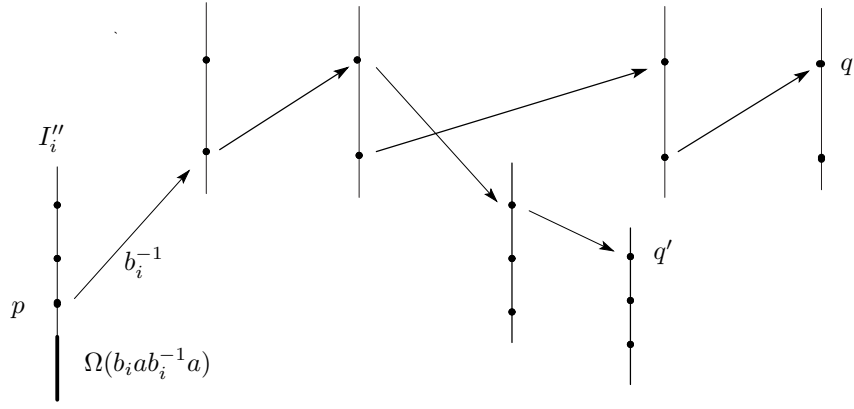

\smallskip\noindent
{\sc Case 4}. $J=I_i'''$, Since 
${\hat W}(\partial_-I_i''')=a{\hat W}(\partial_-I_i''')=ab_i\cdots$,
we have $\partial_-I_i'''\in \Omega(ab_iab_i^{-1})$ and the map in
concern is $b_iab_i^{-1}a$.
This case is similar to Case 3 and is omitted. 
\qed

\bigskip
For the proof of Theorem \ref{t1},
we need Proposition 3.3 of \cite{MR}.
The statement below is adapted for our purpose, and we shall give a complete
proof in Appendix.

\begin{proposition}\label{p61} Let $\phi$ be a dynamical realization of
a circular order $c$ based at $x_0$.
Given any neighbourhood $U$ of $\phi$ in ${\rm Hom}(G,{\rm  Homeo_+(S^1)})$,
there is a neighbourhood $V$ of $c$ in $CO(G)$ such that any order
in $V$ has a dynamical realization based at $x_0$ contained in $U$.
\end{proposition}

\newcommand{\h}{{h_{\rm min}}}
\newcommand{\V}{{\mathcal V}}
\newcommand{\K}{{\mathcal K}}
\newcommand{\E}{{\mathcal E}}

We also need the following proposition whose proof is 
again contained in Appendix. 

\begin{proposition}
 \label{p62}
The space ${\rm Hom}(G,{\rm  Homeo_+(S^1)})$ is locally pathwise
 connected. In fact, more can be said: any neighbourhood of any point contains a pathwise
connected neighbourhood.
\end{proposition}

\bigskip Let us finish the proof of Theorem \ref{t1}.
Recall that  $\phi_\M$ is
a dynamical realization of the circular order
$c_\M$ given by Lemma \ref{l69} based at a point 
$x_0\in I_1$. 
 One can choose $x_0$ in
${\rm Int}(S^1\setminus \cup N_i^{\pm})$, where $N_i^{\pm}$
 is from Lemma \ref{l611}. 
Choose a pathwise connected neighbourhood $U$ of $\phi_\M$ in 
${\rm Hom}(G,{\rm Homeo_+(S^1)})$
 such that
any action $\psi\in U$ satisfies 
$h_i(\psi)(S^1\setminus N_i^{-})\subset N_i^+$ (and eqivalently
$h_i(\psi)^{-1}(S^1\setminus N_i^+)\subset N_i^{-}$),
where
$h_1(\psi)=\psi(\alpha\beta\alpha\beta^{-1})$ and
$h_2(\psi)=\psi(\alpha\beta^{-1}\alpha\beta)$.

\rm
The pingpong lemma (Klein's criterion) asserts that
the subgroup $[G,G]$  generated by $\alpha\beta\alpha\beta^{-1}$
 and $\alpha\beta^{-1}\alpha\beta$
 acts freely at $x_0$ by any $\psi\in U$.
We claim that $G$ itself acts freely at $x_0$ by $\psi$. In fact, if $H$
is the stabilizer at $x_0$, then $H\cap[G,G]=\{e\}$.
That is, the canonical projection $G\to G/[G,G]\cong \Z/6\Z$
is injective on $H$, and hence $H$, if nontrivial, contains an element $\gamma$
of finite order. But $\gamma$ is conjugate either to $\alpha$ or to
$\beta^{\pm1}$, and therefore $\psi(\gamma)$ is fixed point free.
A contradiction shows the claim.

Now let
$V$ be a neighbourhood of $c_\M$ in $CO(G)$ such that a dynamical
realization $\psi$ of an arbitrary order $c'$ of $V$ belongs to $U$.
Then there is a path $\psi_t$, $0\leq t\leq1$, joining
$\phi_\M$ and $\psi$ contained in $U$.
Since all the $\psi_t$ act freely at $x_0$, the circular orders of $G$
obtained from the orbit $\psi_t(G)x_0$ are all the same, independent of $t$.
(Recall that $CO(G)$ is totally disconnected.)
 This shows $c'=c_\M$. The proof of Theorem \ref{t1} is now
complete. \qed

\bigskip
{\sc Proof of Theorem \ref{T}}.
We have constructed a map from the
set of the isolated circular orders to the set of the equivalence classes of
 the Markov systems in Part II, and 
 a map in the reverse
direction in Part III. It is clear from the constructions that one is the inverse
of the other. Thus the proof is done of Theorem \ref{T}. 
\qed

\bigskip
\begin{center}
 \sc\Large Part IV
\end{center}

We present some examples of isolated circular orders of $G$.

\section{Primary examples}

We revisit isolated circular orders of $G$ constructed in \cite{Ma}.

\medskip
\sc Standard example. \rm Here is a Markov system $\M_1$ with
multiplicity 1.
  Just place
intervals $[a],[b],[b^{-1}]$ in the anti-clockwise order on $S^1$.
Clearly there are an involution $a$ of $S^1$ which interchanges  $[a]$ 
and $[[b]]$ and  a period 3 homeomorphism $b$ which circulates
$[a]$, $[b]$ and $[b^{-1}]$. Transitivity of  gaps is
also clear. Thus we obtain an isolated circular order $c_1=c_{\M_1}$.
If the placement of $[a],[b],[b^{-1}]$ is clockwise, we get
another order $c'$. But this is in the same automorphism class of
$c_1$: $c'=(\sigma_0)_*c_1$ for $\sigma_0$ in Proposition \ref{p1} (4).
By the homeomorphism 
$CO(G)\approx LO(B_3)$, $c_1$ corresponds to the Dubrovina-Dubrovin
order \cite{DD}.

\medskip
{\sc Finite lifts of the standard example.}
Let $\phi_1$ be a dynamical realization of the previous example $c_1$.
We shall construct more examples starting with finite lifts of
$\phi_1$.
For $k>1$, let $p_k:S^1\to S^1$ be the $k$-fold covering map.
 A $G$-action $\phi_k$ on $S^1$ is called a $k$-fold lift of $\phi_1$ if
 it satisfies $p_k\phi_k(g)=\phi_1(g)$ for any $g\in G$.
Then the rotation numbers satisfy 
$k\,{\rm rot}(\phi_k(g))={\rm rot}(\phi_1(g))$.
This shows that a lift of an involution is an involution if and only if $k$ is
odd. Likewise a lift of a 3-periodic map is
3-periodic if and only if $k$ is coprime to 3. Therefore
a $k$-fold lift $\phi_k$ of $\phi_1$ exists if and only if $k\equiv
\pm1$ (6):
moreover it is unique if it exists. 
The map $\phi_k(\beta)$ is 3-periodic and has rotation number 
$\pm1/3$ according as $k\equiv\pm1$ (6),
 since 
$(6\ell\pm1)\cdot(\epsilon/3)\equiv 1/3$ (1) implies $\epsilon=\pm1$. 

We shall show that the lift $\phi_k$ is associated with a Markov system
$\M_k$ for $k=6\ell+1$.
The case $k=6\ell-1$ is similar and is left to the reader.
Let $[a]_i$, $[b]_i$ and $[b^{-1}]_i$ ($i\in\Z/k\Z$) be the lifts
of $[a]$, $[b]$
and $[b^{-1}]$ of the privious example $\M_1$ by $p_k$, ordered anticlockwise in $S^1$ as
\begin{equation}
 \label{e71}
\cdots [a]_i,[b]_i,[b^{-1}]_i,[a]_{i+1}\cdots . \ \ \ 
\end{equation}
 The maps $a$, $b$ of the system $\M_k$ is to be
$\phi_k(\alpha)$ and $\phi_k(\beta)$. 
(\ref{e71}) has $3k=18\ell+3$ terms, and $b$, being the lift of $\phi_1(\beta)$,  maps each term to
the term $6\ell+1$ right to it. Therefore we have
\begin{equation}
 \label{e72}
b[a]_i=[b]_{i+2\ell}, \ \ b[b]_i=[b^{-1}]_{i+2\ell}, \ \
b[b^{-1}]_i=[a]_{i+2\ell+1}.
\end{equation}
The sequence (\ref{e71}) is contracted into a sequence
\begin{equation*}
 \label{e73}
\cdots[a]_i,[[b]]_i,[a]_{i+1},[[b]]_{i+1}\cdots
\end{equation*}
of $12\ell+2$ terms. The map $a$ send each term to the term $6\ell+1$
right to it. That is,
\begin{equation}
 \label{e74}
a[a]_i=[[b]]_{i+3\ell},\ \ \ a[[b]]_i=[a]_{i+3\ell+1}.
\end{equation}
 
In order to check (E), we shall classify all the principal gaps into two families $J_i$'s and $J_i'$'s
($i\in\Z/(6\ell+1)\Z$) by the following ordering:
$$
\cdots,[a]_i,J_i,[b]_i,[b^{-1}]_i,J_i',[a]_{i+1},\cdots.$$

By (\ref{e72}) and (\ref{e74}), we get
\begin{equation}
 \label{e75}
b(J_i')=J_{i+2\ell+1}\ \ \ \mbox{ and}
\end{equation}
\begin{equation}
 \label{e76}
ab(J_i')=J'_{i+5\ell+1}.
\end{equation}
Now (\ref{e76}) shows that the group generated by $ab$ acts transitively on the
family $J_i'$'s, since $(6\ell+1,5\ell+1)=1$, and (\ref{e75}) shows that
any $J_i$ from the other family is mapped by $b^{-1}$ to an element of this family. This shows (E).

The circular order defined by $\M_k$ is denoted by $c_k$.
They are from distinct automorphism classes since ${\rm deg}(c_k)=k$.

\section{Further example}

\newcommand{\ua}{{\underline{a}}}
\newcommand{\ub}{{\underline{b}}}
\newcommand{\ubb}{{\underline{b}^{-1}}}

Notice that the ${\rm deg}(c)$ of any isolated circular order $c$ 
is odd, since the involution $a$ transposes $[a]$
and $[[b]]$.
According to our calculation, there is no new
examples of isolated circular orders up to degree $\leq7$. But there is
one in degree 9, which we shall present below.

This example is not well ordered as the
previous one and it is
no use to give an index to each component. Any component of $[a]$
is denoted by the same letter $\ua$. Likewise we use notations $\ub$ and $\ubb$.
Also a component of $[[b]]$ is denoted by $[\ub]$.
Consider the following ordering of 27 intervals in $S^1$ which is
grouped into three:
\begin{equation}
 \label{e81}
\ua\,\ubb\ua\,\ubb\ua\,\ubb\ub\,\ua\,\ubb \mid \ub\,\ua\,\ub\,\ua\,\ub\,\ua\,\ubb\ub\ua \mid
\ubb\ub\,\ubb\ub\,\ubb\ub\,\ua\,\ubb\ub.
\end{equation}
Define a period 3 homeomorphism $b$ of $S^1$ by permuting the three
groups cyclically to the right. Notice that $b$ so defined
 satisfies $b\ua=\ub$, $b\ub=\ubb$ and
$b\ubb=\ua$.
The sequence (\ref{e81}) is contracted to
\begin{equation}
\label{e82}
 \ua\,[\ub]\,\ua\,[\ub]\,\ua\,[\ub]\,\ua\,[\ub]\,\ua\mid[\ub]\,\ua\,[\ub]\,
\ua\,[\ub]\,\ua\,[\ub]\,\ua\,[\ub]
\end{equation}
Define an involution $a$ by transposing the groups.
\newcommand{\ve}[2]
{\left[ \begin{array}{c} #1 \\ #2 \end{array} \right]}
The maps $a$ and $b$ satisfies (A)--(D). To show (E), indicate the
prinicipal gaps by $\ve{i}{j}$ and complementary gaps by $\ve{i}{*}$ as
follows.
$$
\ua\ve{1}{1}\ubb\ve{2}{2}\ua\ve{3}{3}\ubb\ve{4}{4}\ua\ve{5}{5}\ubb\ve{6}{*}\ub\ve{7}{6}\ua\ve{8}{7}\ubb\ve{9}{*}
$$
$$
\ub\ve{1}{8}\ua\ve{2}{9}\ub\ve{3}{1}\ua\ve{4}{2}\ub\ve{5}{3}\ua\ve{6}{4}
\ubb\ve{7}{*}\ub\ve{8}{5}\ua\ve{9}{6}
$$
$$
\ubb\ve{1}{*}\ub\ve{2}{*}\ubb\ve{3}{*}\ub\ve{4}{*}\ubb\ve{5}{*}
\ub\ve{6}{7}\ua\ve{7}{8}
\ubb\ve{8}{*}\ub\ve{9}{9}
$$
Notice that either one of maps $b^{\pm1}$  sends $\ve{i}{j}$ to some 
$\ve{i}{k}$ and the other
to $\ve{i}{*}$, and the map $a$ sends $\ve{i}{j}$
to some $\ve{k}{j}$. We can find a cycle of principal gaps:
$$
\ve{1}{1}\stackrel{a}{\to}\ve{3}{1}\stackrel{b^{-1}}{\to}\ve{3}{3}\stackrel{a}{\to}
\ve{5}{3}\stackrel{b^{-1}}{\to}\ve{5}{5}\stackrel{a}{\to}\ve{8}{5}\stackrel{b^{-1}}{\to}
\ve{8}{7}\stackrel{a}{\to}\ve{6}{7}\stackrel{b^{-1}}{\to}\ve{6}{4}\stackrel{a}{\to}
$$
$$
\ve{4}{4}\stackrel{b}{\to}\ve{4}{2}\stackrel{a}{\to}\ve{2}{2}\stackrel{b}{\to}
\ve{2}{9}\stackrel{a}{\to}\ve{9}{9}\stackrel{b^{-1}}{\to}\ve{9}{6}\stackrel{a}{\to}
\ve{7}{6}\stackrel{b^{-1}}{\to}\ve{7}{8}\stackrel{a}{\to}\ve{1}{8}\stackrel{b^{-1}}{\to}
$$
Therefore (E) is also satisfied. This yields an isolated
circular order $c_9$ of degree 9. It is not in the automorphism classes
of the previous examples. If we arrange the intervals in the clockwise order,
we obtain another order $c'$. But again $c'=(\sigma_0)^*c_9$.

\section{Appendix}

We shall give proofs of Propositions 12.4 and 12.5. The former holds
true for an arbitrary countable group $\oG$.

\bigskip
{\sc Proposition 12.4}. \em Let $\phi$ be a dynamical realization of
a circular order $c$ based at $x_0$.
Given any neighbourhood $U$ of $\phi$ in ${\rm Hom}(\oG,{\rm  Homeo_+(S^1)})$,
there is a neighbourhood $V$ of $c$ in $CO(\oG)$ such that any order
in $V$ has a dynamical realization based at $x_0$ contained in $U$.
\rm \bigskip

\bd  Given a finite set
$F$ of $\oG$ and $\epsilon>0$, let
$$
U(F,\epsilon,\phi)=\{\psi\in{\rm Hom}(\oG,{\rm  Homeo_+(S^1)})\mid
\Abs{\psi(g)-\phi(g)}_0<\epsilon,\ \forall g\in F\}.$$
One may replace the given $U$ in the proposition by a smaller neighbourhood
$U(F,\epsilon,\phi)$. 
One may also assume that $F$ is symmetric: $g\in F$ implies $g^{-1}\in F$.
Given a finite subset $S$ of $\oG$,
define a neighbourhood $V_S(c)$ of $c$  in $CO(H)$ by
$$
V_S(c)=\{c'\in CO(\oG)\mid
c'\vert_{(S\cup\{e\})^3}=c\vert_{(S\cup\{e\})^3}\}.
$$
Then any $c'\in V_S(c)$ admits a dynamical realization $\psi_{c',S}$
based at $x_0$
such that $\psi_{c',S}(g)x_0=\phi(g)x_0$ for any $g\in S$.
Our aim is first to find good $S$ and then to alter $\psi_{c',S}$ by a $x_0$-preserving
conjugacy to obtain a homomorphism contained in $U(F,\epsilon,\phi)$.
 Choose $\delta>0$ so that
whenever $f\in F$ and  $\abs{x-y}<\delta$, we have
$\abs{\phi(f)x-\phi(f)y}<\epsilon$.

\smallskip
{\sc Case 1}, \em  $\phi(\oG)x_0$ is dense in $S^1$.
\rm
Choose a finite subset $S'$
 of $\oG$ 
so that $X=\phi(S')x_0$ is $\delta/2$-dense in $S'$ and  define
$S=FS'\cup S'$. In this case there is no need for the alteration:
we shall show that 
$\psi_{c',S}\in  U(F,\epsilon,\phi)$. 
In fact, for any $x\in X$, we have
 $x=\psi_{c',S}(g)x_0=\phi(g)x_0$ for $g\in S'$. Moreover $\psi_{c',S}(f)x=\phi(f)x$ for
any $x\in X$ and $f\in F$, because $\phi(f)x=\phi(fg)x_0$, $\psi_{c',S}(f)x=\psi_{c',S}(fg)x_0$
 and $fg\in S$. Choose any $f\in F$ and any $y\in S^1$.
Then $y$ belongs to the closure $\Cl(J)$ of some gap
$J$ of $X$.
Since $X$ is $\delta/2$-dense, we have $\abs{J}<\delta$,
and hence $\abs{\phi(f)J}<\epsilon$. Since both points $\phi(f)y$ and
$\psi_{c',S}(f)y$ belong to the same interval $\phi(f)\Cl(J)=\psi_{c',S}(f)\Cl(J)$, we have 
$\abs{\phi(f)y-\psi_{c',S}(f)y}<\epsilon$. Since $f\in F$ and $y\in S^1$ are
arbitrary, we obtain 
$\psi_{c',S}\in U(F,\epsilon,\phi)$, as is required.

\smallskip
{\sc Case 2}. \em $\phi(\oG)x_0$ is not dense in $S^1$.
\rm
Denote by $J_0$ the gap of $\Cl(\phi(\oG)x_0)$ whose left end point is $x_0$.
Define $\h\in \oG$ so that $\phi(\h)x_0$ is the right end point of $J_0$.
Since $\phi$ is tight at $x_0$,  any gap $J$ of $\Cl(\phi(\oG)x_0)$
is a translate of $J_0$.
Denote by $\mathcal V$ the set of
 the gaps $J$  such that $\abs{J}\geq\delta$. Let
$$S_1=\{g\in\oG\mid \phi(g)(J_0)\in\mathcal V\}$$
and let $S_2=S_1\cup S_1\h$.
Thus the end points of any interval of $\V$ belongs to $\phi(S_2)x_0$.
Add some more elements to $S_2$ to form a finite subset $S_3$ such that 
$\phi(S_3)x_0$ is $\delta/2$-dense in $\Cl(\phi(\oG)x_0)$.
 Notice that any gap of $\phi(S_3)x_0$ either belongs to $\V$
or is of length $<\delta$.
Finally let $S=FS_3\cup S_3$.
We have the following property by the same argument as in Case 1.

\smallskip
\em (1)  For any point  $x\in\phi(S_3)x_0$ 
and $f\in F$, we have $\psi_{c',S}(f)x=\phi(f)x$.\rm

\smallskip\noindent
Define a one dimensional simplicial complex $\K$ as follows. The
vertex set of $\K$ is $\V$. Two vertices $J$ and $J'$ are joined by an edge
if there is $f\in F$ such that $\phi(f)J=J'$. 
Such $f$ is unique since $\phi(\oG)$ acts simply transitively on the gaps
of  $\Cl(\phi(\oG)x_0)$: $\phi(g)J_0=J_0$ implies $g=e$.
We label the direted edge from $J$ to $J'$ by $f$.
For a directed ege $J\stackrel{f}{\to}J'$ of $\K$, we have
$\psi_{c',S}(f)(J)=J'$ by (1).
The simple transitivity on the gaps shows the following.

\smallskip \em
(2) For a directed cycle
$$J_1\stackrel{f_1}{\to}J_2\stackrel{f_2}{\to}\cdots
J_{n}\stackrel{f_n}{\to}J_{n+1}=J_1,
$$
we have $f_n\cdots f_2f_1=e$.

\smallskip\rm

Let us consider $h\in{\rm Homeo}_+(S^1)$ supported on the union 
of the closures of the intervals in $\V$ and leaving the end points of
the intervals fixed.
We claim that there is $h$ such 
that for any directed edge $J\stackrel{f}{\to}J'$
of $\K$, $h\psi_{c',S}(f)h^{-1}=\phi(f)$ on $J$. To show this, consider a
spanning tree $T_\nu$ of each component $\K_\nu$ of $\K$.
 Define $h$ to be the identity
on a prescribed base vertex $J_1$ of $T_\nu$. For any directed
edge  $J_1\stackrel{f}{\to}J_2$ of $T_\nu$, define $h$ on $J_2$ so that
$h\psi_{c',S}(f)h^{-1}=\phi(f)$ holds on $J_1$. We continue this process
along directed paths in $T_\nu$ issuing at $J_1$ until we define $h$ on
all the vertices of $\K_\nu$. Then for any
directed edge $J\stackrel{f}{\to}J'$
of $T_\nu$, $h\psi_{c',S}(f)h^{-1}=\phi(f)$ on $J$ (even if the edge is
directed toward the base vertex).
Also for an edge of $\K_\nu$ not in $T_\nu$, we have the same
equality thanks to the relation (2). The proof of the claim is over.
 
Let us define a homomorphsim $\psi$ by  $\psi(g)=h\psi_{c',S}(g)h^{-1}$
 for
any $g\in\oG$. 
The homomorphism $\psi$ still satisfies (1) with $\psi_{c',S}$ replaced
 by $\psi$.
Let $J$ be any gap of 
$\phi(S_3)x_0=\psi(S_3)x_0$ and $f$ any element of $F$.
If either $\abs{J}<\delta$ or $\abs{\phi(f)J}<\delta$, then
 $\psi(f)$ is $\epsilon$-near to $\phi(f)$ on $J$.
If not, $J\stackrel{f}{\to}\phi(f)J$ is an edge of
$\K$, and $\psi(f)=\phi(f)$ on $J$. This shows $\psi\in
 U(F,\epsilon,\phi)$, as is required.
 \qed

\bigskip
{\sc Proposition 12.5}. \em
The space ${\rm Hom}(G,{\rm  Homeo_+(S^1)})$ is locally pathwise
 connected. In fact, any neighbourhood of any point contains a pathwise
connected neighbourhood.
\rm \bigskip

\bd The space in question is homeomorphic to $Q_2\times Q_3$, where
\begin{equation*}
\begin{array}{c}
Q_2=\{a\in{\rm  Homeo}_+(S^1)\mid a^2=id\}\\
Q_3=\{b\in{\rm  Homeo}_+(S^1)\mid b^3=id\}.
\end{array}
\end{equation*}
 Choose a base point $x_0\in S^1$.
A map $a$ in $Q_2$ is specified by a point $a(x_0)$ and an orientation
preserving homeomorphism from $[x_0,a(x_0)]$ to $[a(x_0),x_0]$. So $Q_2$ is
homeomophic to $(0,1)\times{\rm  Homeo}_+(0,1)$.
Likewise $Q_3$ is homeomorphic to 
$(0,1)\times{\rm  Homeo}_+(0,1)\times{\rm  Homeo}_+(0,1)$.
It suffices to show that the space ${\rm  Homeo}_+(0,1)$ satisfies
the claimed property. For any neighbourhood $U$ of a point 
$f\in{\rm  Homeo}_+(0,1)$, choose $\epsilon>0$ so that
$$
U(f,\epsilon)=\{g\in{\rm  Homeo}_+(0,1)\mid \Abs{g-f}<\epsilon\}$$
is contained in $U$. Then 
for any 
$g\in U(f,\epsilon)$,
the path $\{(1-t)f+tg\}_{0\leq t\leq1}$ is
contained in $U(f,\epsilon)$. \qed

\end{document}